\documentclass[reqno]{amsproc}

\usepackage{graphicx}

\usepackage{hyperref}

\hypersetup{
colorlinks=true,
linkcolor=blue,
anchorcolor=blue,
citecolor=blue
}

\usepackage{amssymb}
\usepackage{amsfonts}
\usepackage{amsmath}

\newtheorem{theorem}{Theorem}[section]
\newtheorem{lemma}[theorem]{Lemma}
\newtheorem{proposition}[theorem]{Proposition}

\theoremstyle{definition}
\newtheorem{definition}[theorem]{Definition}

\newtheorem{remark}[theorem]{Remark}

\numberwithin{equation}{section}

\begin{document}

\title[]{John-Nirenberg inequalities for noncommutative column BMO and Lipschitz martingales}

\author{Guixiang Hong}
\address{Guixiang Hong, Institute for Advanced Study in Mathematics, Harbin Institute of Technology, Harbin 150001, China}
\email{gxhong@hit.edu.cn}
\thanks{The first and third authors were supported in part by NSFC No.12071355.}

\author{Congbian Ma}
\address{School of Mathematics and Statistics, Xinxiang University, Henan, China}
\email{congbianm@whu.edu.cn}
\thanks{The second author was supported in part by NSFC No.11801489.}

\author{Yu Wang}
\address{School of Mathematics and Statistics, Wuhan University, Wuhan, China}
\email{yuwang1995@whu.edu.cn}





\begin{abstract}
In this paper, we continue the study of John-Nirenberg theorems for BMO/Lipschitz spaces in the noncommutative martingale setting. As conjectured from the classical case, a desired noncommutative ``stopping time" argument was discovered to obtain the distribution function inequality form of John-Nirenberg theorem. This not only provides another approach without using duality and interpolation to the results for spaces $\mathsf{bmo}^c(\mathcal M)$ and ${{\Lambda}^{{c}}_{\beta}}(\mathcal{M})$, but also allows us to find the desired version of John-Nirenberg inequalities for
spaces $\mathcal{BMO}^c(\mathcal M)$ and ${\mathcal L^{{c}}_{\beta}}(\mathcal{M})$. And thus we solve two open questions after \cite{ref5, ref3}. As an application, we show that Lipschitz space is also the dual space of noncommutative Hardy space defined via symmetric atoms. Finally, our results for ${\mathcal L^{{c}}_{\beta}}(\mathcal{M})$ as well as the approach seem new even going back to the classical setting.
\end{abstract}

\maketitle

\section{Introduction}

The purpose of this paper is to study the John-Nirenberg {theorems} of noncommutative column BMO and Lipschitz martingales. In the seminal paper \cite{ref35}, Pisier and Xu introduced the noncommutative martingale framework and established the noncommutative Burkholder-Gundy inequalities and Fefferman-Stein duality theorem between Hardy and BMO spaces. Later on, Junge  \cite{ref30} made a breakthough on noncommutative Doob's maximal inequality. Based on the two fundamental papers, the theory of noncommutative martingales has been rapidly developed. We refer the reader to \cite{ref32, ref7} for noncommutative Burkholder/Rosenthal inequalities and ergodic theorems, to \cite{ref36, ref38} for the noncommutative Gundy and Davis decompositions, to \cite{ref13,  ref33} for noncommutative martingale transform and conditioned square functions, to \cite{JRWZ20, JOW18, ref6} and the references therein for noncommutative differential subordinates and good-$\lambda$ inequalities. We should also mention another two works directly related with the objectives of this paper, which are
the  John-Nirenberg theorems  for noncommutative martingale BMO spaces by Mei and the first author \cite{ref5} (see also \cite{ref19,ref18}) and the atomic decomposition for Hardy spaces with small exponents by Chen,
Randrianantoanina and Xu \cite{ref3} (see also \cite{ref2}).

Let $\mathcal M$ be a von Neumann algebra  with a
normal faithful tracial state $\tau$.  Let
$(\mathcal M_n)_{n\geq1}$ be an increasing sequence of von Neumann
subalgebras of $\mathcal M$ such that the union of $\mathcal M_n$'s is $w^*$-dense
in $\mathcal M$ and $\mathcal{E}_n$'s be the associated conditional expectations. A sequence $(x_n)_{n\geq1}$ is called a $L_p$-martingale if $x_n\in L_p(\mathcal M)$ and $\mathcal E_n(x_{n+1})=x_n$ for each $n\geq1$. The column space $\mathcal{BMO}^c(\mathcal M)$ is defined to be subset of $L_2(\mathcal M)$ (or equivalently $L_2$-martingales) with
$$\|x\|_{\mathcal{BMO}^c}:=\sup_{n\geq1}
\|\mathcal{E}_n|x-\mathcal E_{n-1}(x)|^2\|^{\frac 12}_{\infty}<\infty;$$
The mixture space is then defined as
 $$\mathcal{BMO}(\mathcal M)=\{x\in L_2(\mathcal M): \|x\|_{\mathcal{BMO}}<\infty\}$$
with
$$\|x\|_{\mathcal{BMO}}=\max\{\|x\|_{\mathcal{BMO}^c},\|x^*\|_{\mathcal{BMO}^c}\}.$$
Keeping in mind that the classical method based on the distribution function inequalities {is} very difficult to be adapted to the noncommutative setting due to the noncommutativity of operator product and the lack of an efficient analogue of the notion of stopping times, motivated by  \cite{ref18}, Mei and the first author in \cite{ref5} proved the John-Nirenberg theorem for $\mathcal {BMO}(\mathcal M)$ starting with the following form: for $0<p<\infty$, we have
\begin{align}\label{JNforBMO}\alpha^{-1}_p\|x\|_{\mathcal{BMO}}\leq\mathcal{PB}_p(x)\leq\beta_p\|x\|_{\mathcal{BMO}},\end{align}
where
\begin{align}\label{PB}
\mathcal{PB}_p(x)=\sup_n\sup_{e\in \mathcal{P}(\mathcal M_n)} \max\{
\|(x-x_{n-1})\frac{e}{(\tau(e))^{1/p}}\|_{p},
\|\frac{e}{(\tau(e))^{1/p}}(x-x_{n-1})\|_{p}\}.
\end{align}
Here $\mathcal P(\mathcal M)$ denotes the lattice of projections of $\mathcal M_n$ and
the two constants $\alpha_p$ and $\beta_p$ have the following
properties
\begin{center} $\quad\quad\quad\alpha_p=1$ for $2\leq p<\infty$, \ \ $\alpha_p\leq C^{1/p-1/2}$ for $0< p<2$,\end{center}
\begin{center} $\beta_p\leq cp$ for $2\leq p<\infty$, \ \  $\beta_p=1$ for $0< p<2$,\end{center}
where $c$ and $C$ are two universal positive constants.

The above form \eqref{JNforBMO}, also called a {\it fine} version in \cite{ref5} which corresponds perfectly to the classical case, implies in a standard way another two forms of John-Nirenberg theorem in terms of the distribution inequality and the exponential integrability (see e.g. \cite{ref5} for the details). However regarding the column space $\mathcal{BMO}^c(\mathcal M)$, the situation is much more complicated. Indeed, an exact version in terms of $L_p$ spaces corresponding to \eqref{JNforBMO} was disproved in \cite[Remark 3.14]{ref5} and there holds only for $2\leq p<\infty$ a {\it crude} version in terms of Hardy spaces $\mathcal H^c_p$ (cf \cite[Theorem 3.8]{ref5}), that is, projections $e$'s and $L_p$-norm in \eqref{PB} are replaced by elements $b\in L_p$'s and $\mathcal H^c_p$-norm respectively. The reason behind is that the restricted scale $p\geq2$ does not allow us to exploit an extreme point property of $L_q$ spaces for $q\leq1$ to modify the crude version to the fine version, see the proof of \cite[Theroem 3.16]{ref5} for the details of the argument. So the question on the John-Nirenberg theorem for $\mathcal{BMO}^c(\mathcal M)$ was open before the present paper.

In the same paper \cite{ref5}, based on a similar idea depending on BMO interpolation, the authors established also John-Nirenberg theorem for the conditioned space $\mathsf{bmo}^c(\mathcal M)$. More recently, starting with the atomic decomposition of Hardy spaces, Chen, Randrianantoanina and Xu \cite{ref3}, via the Fefferman-Stein duality, obtained the fine form of John-Nirenberg inequality for the conditioned Lipschitz space $\Lambda^c_\beta(\mathcal M)$ defined as
\begin{align*}
{\Lambda}^c_{\beta}(\mathcal{M})=&\{x\in L_2(\mathcal{M}):\;\|x\|_{{\Lambda}^c_{\beta}}<\infty\}
\end{align*}
with
$$
\|x\|_{{\Lambda}^c_{\beta}}:=\max\{\|\mathcal{E}_1(x)\|_\infty, \;\sup_{n\geq1}\sup_{e\in{\mathcal{P}(\mathcal{M}_n)}}{\frac{\|(x-x_n)e\|_{\mathsf{h}_p^c}}{(\tau(e))^{\beta}}}\}$$
which is exactly $\mathsf{bmo}^c(\mathcal M)$ when $\beta=0$.  We refer the reader to Theorem \ref{thm:JNlambda} (iii) for the above-mentioned fine version of $\Lambda^c_\beta(\mathcal M)$, and to the body of the paper for the notation or notion not defined here or below. However, inspired by the classical case, before the present paper it had been an open problem to have a direct or constructive proof of the noncommutative John-Nirenberg theorem without involving BMO interpolation or Fefferman-Stein duality, which should have more applications in noncommutative analysis.

\smallskip
In the present paper, we solve the above two open problems. More precisely, we establish the fine version of John-Nirenberg theorem for $\mathcal{BMO}^c(\mathcal M)$ which corresponds perfectly to the classical result (see e.g. \cite[Theorem 2.50]{ref25}): let $x\in\mathcal{BMO}^c(\mathcal{M})$ and $0< p<\infty$, then
\begin{align}\label{JNfp}
\alpha_p^{-1}\|x\|_{\mathcal{BMO}^c}\leq\|x\|_{\mathcal{BMO}^c_p}\leq\beta_p\|x\|_{\mathcal{BMO}^c},
\end{align}
where
$$
\|x\|_{\mathcal{BMO}^c_{p}}= \sup_{n\geq1}\sup_{e\in{\mathcal{P}(\mathcal{M}_n)}}{\frac{\|S_c(x-x_{n-1})e\|_{p}}{(\tau(e))^{\frac{1}{p}}}}
$$
and $S_c(y)$ is the column square function defined as
$$S_c(y)=\Big(\sum^\infty_{k=1}|dy_k|^2\Big)^\frac12$$
with $dy_k=\mathcal E_k(y)-\mathcal E_{k-1}(y)$ for $k\geq2$ and $dy_1=\mathcal E_1(y)$;
moreover, our approach to \eqref{JNfp} is via the distribution function inequality form:
 let $x\in \mathcal{BMO}^c(\mathcal M)$, then for any $n\geq1$, $P\in\mathcal P(\mathcal M_n)$ and $\lambda\geq0$, there holds
\begin{align}\label{JNfd}
\tau \left(I_{[\lambda,\infty)}\big((P S_c^2 (x-x_{n-1}) P)^{\frac{1}{2}} \big)\right)  \le 2(1-e^{-2})^{-1} e^2 e^{- \frac{ \lambda}{e{\| x \| }_{{{\mathcal{BMO}^c}}}}} \tau(P),
\end{align}
which in turn follows from a noncommutative ``stopping time" argument that has been highly desired but unreachable before. Here, ${I}_{[\lambda,\infty)}{(a)}$ denotes the spectral projection of $a\in \mathcal M$ corresponding to the interval $[\lambda,\infty)$.

Moreover, this approach is strong and flexible enough so that it not only provides another method to deal with John-Nirenberg theorem for the conditioned space $\Lambda^c_\beta(\mathcal M)$, but also applies to the Lipschitz space $\mathcal{L}^c_{\beta}(\mathcal{M})$ defined as
$$\mathcal{L}^c_{\beta}(\mathcal{M})=\{x\in L_2(\mathcal{M}): \|x\|_{\mathcal{L}^c_{\beta}}<\infty\}$$
with
$$
\|x\|_{\mathcal{L}^c_{\beta}}:= \sup_{n\geq1}\sup_{e\in{\mathcal{P}(\mathcal{M}_n)}}{\frac{\|(x-x_{n-1})e\|_{2}}{(\tau(e))^{\beta+\frac{1}{2}}}},
$$
which is exactly $\mathcal {BMO}^c(\mathcal M)$ when $\beta=0$. For this reason, we will state and prove the John-Nirenberg theorem directly for  $\mathcal{L}^c_{\beta}(\mathcal{M})$ (see Section \ref{section 2}) and $\Lambda^c_\beta(\mathcal M)$ (see Theorem \ref{thm:JNlambda}). To the authors' best knowledge, these results are new even in the case $\mathcal M$ being commutative and $\beta>0$. As a further application, in Section \ref{section 3},  we  obtain the moment characterization of $\Lambda^c_\beta(\mathcal M)$ in terms of symmetric spaces in a more direct way which is much clearer than the interpolation arguments given in \cite{ref19}, and then show in Section \ref{section 4} that the noncommutative Hardy space ${\mathsf{h}}^c_{p,E} (\mathcal M)$ ($0 < p\leq 1)$ defined via symmetric space atoms is also a predual space of Lipshitz space $\Lambda^c_{1/p-1}(\mathcal M)$.

\smallskip

To explain the difficulties and gain some intuition about our approach, let us recall below briefly the arguments for classical BMO space on a probability space $(\Omega,\mathcal F,\mathbb P)$ equipped with sigma algebras $(\mathcal F_n)_{n\geq1}$ with associated conditional expectations $(\mathbb E_n)_{n\geq1}$ (see e.g. \cite[Section 2.4]{ref25}). Let $f\in {{\mathcal{BMO}}}(\Omega)$ and $f_n:=\mathbb E_n(f)$. The key estimate is the following: for any  $n\in \mathbb{N}$, any $E\in \mathcal{F}_n$ and any $\lambda,\mu>0$,
\begin{align}\label{distrbution function method}
 \mathbb{P}&(\{\omega\in E: S(f-f_{n-1})\geq\lambda+\mu\}) \notag\\
\leq& \frac{\|f\|_{{{}{\mathcal{BMO}}}}}{\mu}\mathbb{P}(\{\omega\in E: S(f-f_{n-1})\geq\lambda\}),
\end{align}
where
$$S(g)=\left(\sum_{k=1}^\infty|dg_k|^2\right)^{\frac12}.$$
The classical proof of (\ref{distrbution function method}) rests on the construction of two stopping times with respect to $\lambda+\mu$ and $\mu$ for the non-decreasing adapted sequence $(S(f_k-f_{n-1}))_{k\geq n}$  which is obviously absent in the noncommutative setting. Motivated by Cuculescu's construction (see the context after \eqref{Cuc}), one may rewrite its proof as follows which might be
transferable to the noncommutative setting,
\begin{align}\label{distribution functin analysis}
 &\quad\mathbb{P}(\{\omega\in E: S(f-f_{n-1})\geq\lambda+\mu\})\notag\\
&\leq\int_E\big(\chi_{{[\lambda+\mu,\infty)}}(S(f-f_{n-1}))\big)
(\frac{S(f-f_{n-1})-\lambda}{\mu})\notag\\
&\leq\frac1{\mu}\int_E
\big(\chi_{{[\lambda,\infty)}}(S(f-f_{n-1}))\big)(S(f-f_{n-1})-\lambda)\notag\\
&= \frac{1}{\mu}\sum_{{k=n}}^\infty\Big(\int_E
\big(\chi_{{[\lambda,\infty)}}(S(f_k-f_{n-1}))-\chi_{{[\lambda,\infty)}}(S(f_{k-1}-f_{n-1}))\big)\notag\\
&\quad\cdot
(S(f-f_{n-1})-S(f_{k-1}-f_{n-1}))+(S(f_{k-1}-f_{n-1})-\lambda)\Big)\notag\\
&\leq \frac{1}{\mu}\sum_{{k=n}}^\infty\Big(\int_E
\big(\chi_{{[\lambda,\infty)}}(S(f_k-f_{n-1}))-\chi_{{[\lambda,\infty)}}(S(f_{k-1}-f_{n-1}))\big)\notag\\
&\quad\cdot
(S(f-f_{n-1})-S(f_{k-1}-f_{n-1}))\Big)\notag\\
&\leq\frac{1}{\mu}\sum_{{k=n}}^\infty\int_E
\big(\chi_{{[\lambda,\infty)}}(S(f_k-f_{n-1}))-\chi_{{[\lambda,\infty)}}(S(f_{k-1}-f_{n-1}))\big)\notag\\
&\quad\quad\cdot
\|\mathbb E_k(S(f-f_{n-1})-S(f_{k-1}-f_{n-1}))\|_\infty\notag\\
&\leq \frac{1}{\mu}\|f\|_{{{\mathcal{BMO}}}}\sum_{{k=n}}^\infty\int_E
\big(\chi_{{[\lambda,\infty)}}(S(f_k-f_{n-1}))-\chi_{{[\lambda,\infty)}}(S(f_{k-1}-f_{n-1}))\big)\notag\\
&=\frac{1}{\mu}{\|f\|_{{{}{\mathcal{BMO}}}}}\mathbb{P}(\{\omega\in E: S(f-f_{n-1})\geq\lambda\}).
\end{align}
 {Then for a given $\lambda>0$, choosing the nonnegative number  $k$ such that $ke\leq \lambda< (k+1)e$ and applying (\ref{distrbution function method}) repeatedly with $\lambda=(k-1)e,...,e$ and $\mu=e$, one gets the distributional function inequality form of John-Nirenberg theorem for ${{}{\mathcal{BMO}}}(\Omega)$,
\begin{align}\label{JNcommu}
&\quad\mathbb{P}(\{\omega\in E:S(f-f_{n-1})\geq\lambda\})\notag\\
&\leq\mathbb{P}(\{\omega\in E:S(f-f_{n-1})\geq ke\})\notag\\
&\leq (e^{-1})^{{k-1}}\mathbb{P}(E)\leq e^2e^{-\frac{\lambda}{e}}\mathbb{P}(E).
\end{align}}
This give the result of the case $\|f\|_{{{}{\mathcal{BMO}}}}=1$ and the general case follows by scaling.

However, the above ``easy" argument does not admit a noncommutative counterpart due to the failure of some elementary inequalities, and one requires several genuinely new ideas. Let us explain two main ones.

In the first identity of \eqref{distribution functin analysis}, we have used the obvious fact that for any sequence of non-decreasing non-negative functions $(g_k)_{k\geq n}$ with the limit $g_\infty$, the level set $\{g_\infty\geq\lambda\}$ can be expressed as a disjoint union
$$\{g_\infty\geq\lambda\}=\bigcup^\infty_{k= n+1}\{g_k\geq\lambda,\; g_{k-1}<\lambda\}\bigcup \{g_n\geq\lambda\}.$$
The above trivial fact, nevertheless, does not hold for sequences of operators since the characteristic function is not operator monotone. This failure leads to a huge amount of additional work to get a noncommutative analogue of the second inequality of \eqref{distribution functin analysis}, which relates two sequences of projections with respect to two levels. For instance, in order to exploit Cuculescu's construction for martingales, we embed $\mathcal{M}$ into a larger matrix algebra to linearize the sequence of square functions $(S(f_k-f_{n-1}))_{k\geq n}$, and then it is quite technical to deal with the simultaneous appearance of two projections that are not commuting in one expression. We refer to Lemma \ref{lem3.2} and its proof for more discussion and details.

The embedding of $\mathcal M$ into a larger matrix algebra induces an additional difficulty to conclude noncommutative version of \eqref{JNcommu} from an intermediate estimate like \eqref{distrbution function method} since the identity in the matrix algebra is uncontrollable. For this purpose, we need Lemma \ref{lem3.1} to obtain $e_{n+1,n+1}+e_{n+3,n+3}$ that has finite trace instead of the identity on the right hand of \eqref{pro1}. The last but not the least, we need also to explore the extreme point property of $L_1$ to obtain a noncommutative analogue of the last inequality of \eqref{distribution functin analysis}. See the proof of {Proposition \ref{lem3.3}}  and Theorem \ref{lem3.4} for more details.


\section{John-Nirenberg theorem for $\mathcal {L}_\beta^c(\mathcal M)$ and $\Lambda^c_\beta(\mathcal M)$}
\label{section 2}
The main purpose of this section is devoted to the establishment of John-Nirenberg theorem in terms of a distribution function inequality. Then  from that, we deduce easily the exponential integrability and $p$-moment characterization.

For the given noncommutative measure space $(\mathcal{M},\tau)$, denote by $L_0(\mathcal{M})$ the set of all the $\tau$-measurable operators.
 For each $0<p<\infty$, let $L_p(\mathcal{M})$ be the subspace of $x\in L_0(\mathcal{M})$ such that
$$\|x\|_p:=\big(\tau(|x|^p)\big)^{\frac1p}<\infty.$$
When $p=\infty$, we take the operator norm with convention.
For $x\in L_0\mathcal{(M})$, the distribution function $\lambda(x)$ of $x$ is defined by $\lambda_t(x)=\tau(I_{(t,\infty)}(|x|))$ for $t>0$,
and the generalized singular numbers $\mu(x)$ by $\mu_t(x)=\inf\{s>0: \lambda_s(x)\leq t\}$ for $t>0$.

\subsection{Exponential decay of the distribution function for $\mathcal {L}_\beta^c(\mathcal M)$}
\label{subsection 2.1}
 As explained in the introduction, one may encounter lots of difficulties to adapt classical stopping time argument into the noncommutative setting. The approach that we provide below includes several new ideas, and is actually new even going back to the classical setting.

The first idea of our approach is to linearize the square function in a nice way. For this purpose, we will embed $\mathcal{M}$ into a larger von Neumann algebra.
Fix $n\geq1$ and $N> n$. Let $(\Omega,\mathcal{F},\mathbb{P})$ be a probability space and let $\varepsilon_n,\varepsilon_{n+1},...,\varepsilon_N$ be a sequence of independent Rademacher variables. For a fixed $P \in {\mathcal{P}(\mathcal{M}_n)}$, consider the amplification algebra
$\mathcal{N}=\mathbb{M}_{N+2-n} \otimes L_{\infty} (\Omega,\mathcal{F},\mathbb{P})\otimes P\mathcal{M}P$, where $\mathbb{M}_{N+2-n} $ is the algebra of $(N+2-n)\times(N+2-n)$ matrices with the usual trace.   We equip $\mathcal{N}$ with the usual tensor trace $\nu$ and the filtration ${(\mathcal{N}_{m})_{m=n}^N}$ with $\mathcal{N}_{m}={\mathbb{M}_{N+2-n} \otimes L_{\infty} (\Omega,\mathcal{F}_{m},\mathbb{P})\otimes P{\mathcal{M}}_m P}$,
where  $\mathcal{F}_m$ stands for the $\sigma$-algebra generated by the variables $\varepsilon_n,\varepsilon_{n+1},...,\varepsilon_m$. The conditional expectation associated with ${\mathcal{N}}_m $ is denoted by $ \bar{\mathcal{E}}_m$.

For $ x\in \mathcal{L}^{c} _{ \beta }(\mathcal{M})  $ with  $ \beta \ge 0$, we consider the associated sequence $y=(y_m)_{m=n}^N$ given by
\begin{eqnarray}\label{sss}
 y_{m} = \sum_{k=n}^m(e_{n+1,k+2} + e_{k+2,n+1 })  \otimes \varepsilon_{k} \otimes \frac{(P|dx_{k}|^2 P)^\frac{1}{2}} {{\tau(P)}^\beta}, \end{eqnarray}
where $e_{i,j}$ are the standard units of $\mathbb{M}_{N+2-n}$. The associated difference is given by the formula  $dy_n=y_n$ and $dy_m=y_m-y_{m-1}$ for $n<m\leq N$.
Clearly, $y$ is a martingale with respect to $({\mathcal{N}_m})_{m=n}^N.$  For $\lambda \geq 0$ we define the required  sequence $R^\lambda=(R_m^\lambda)_{m\ge n} $ of projections associated with $y$, given by ${R} _{n-1} ^\lambda = \bar{I}= (\begin{matrix}
\sum_{k=n+1}^{N+2} e_{k,k} \end{matrix} )\otimes 1 \otimes P $ and, inductively,
$$ R_m^{\lambda}= R_{m-1}^{\lambda} {\bar{I}}_{(-\infty,\lambda)}(R_{m-1}^{\lambda} y_m R_{m-1}^{\lambda}).$$
Here, $\bar{I}_{(-\infty,\lambda)}{(a)}$ denotes the spectral projection of $a\in L_0(\mathcal N)$ corresponding to the interval $(-\infty,\lambda)$.
It is well-known that these projections enjoy the following properties (see e.g. \cite{ref39}):
\begin{itemize}\label{Cuc}
\item $R_m^{\lambda} \in {\mathcal{N}}_m,\;\forall n\leq m\leq N;$
\item $R_m^{\lambda}$ commutes with $R_{m-1}^{\lambda} y_mR_{m-1}^{\lambda},\;\forall n\leq m\leq N;$
\item $ R_m^{\lambda} y_m R_m^{\lambda} \le  \lambda {{}R_m^{\lambda}},\;\forall n\leq m\leq N.$
\end{itemize}

The following intermediate estimate, as a noncommutative analogue of \eqref{distrbution function method}, will play a key role in establishing the distribution function inequality form of John-Nirenberg theorem.

\begin{proposition}\label{lem3.3}  For any $\lambda, \mu >0$,   one has
\begin{align}\label{pro1} \nu (\bar{I}-R_N^{\lambda+\mu}) \le \frac{4}{{\mu}^2}  \left \| x  \right \|^2_{\mathcal{L}^c_{\beta}} \nu \big((\bar{I}-R_N^{\lambda})((e_{n+1,n+1}+e_{n+3,n+3}) \otimes 1 \otimes P)\big).
\end{align}
\end{proposition}

To show Proposition \ref{lem3.3}, we need two lemmas. The first one is an inequality involving two sequences of projections with different levels, which correspond to two classical stopping times. In the classical setting, this inequality follows easily from commutativity and the Chebyshev inequality. But it requires considerably careful analysis in the noncommutative setting.
\begin{lemma}\label{lem3.2} For any $\lambda, \mu >0$,   we have
$$  \nu (\bar{I}-R_N^{\lambda+\mu}) \le \frac{2}{{\mu}^2} \nu ((\bar{I}-R_N^{\lambda}){(y_N- \lambda \bar{I})}^2).$$
\end{lemma}
The above result has been essentially established in \cite{JOW18, ref6}. But we have a new observation which enables us to simplify greatly the arguments.
\begin{proof}  Fix  $m\in \{n,n+1,...,N\}$. By the definition of $R_m^{\lambda+\mu}$ and using the martingale property of  $y$,
\begin{align*}
&\nu(R_{m-1}^{\lambda+\mu}-R_{m}^{\lambda+\mu})& \\
&=\nu\Big({\bar{I}}_{[\lambda+\mu,\infty)}\big((R_{m-1}^{\lambda+\mu}-R_{m}^{\lambda+\mu})y_m(R_{m-1}^{\lambda+\mu}-R_{m}^{\lambda+\mu})
\big)\Big) \\
&=\nu\Big({\bar{I}}_{[\lambda+\mu,\infty)}\big((R_{m-1}^{\lambda+\mu}-R_{m}^{\lambda+\mu})\bar{\mathcal{E}}_m(y_N)(R_{m-1}^{\lambda+\mu}-
R_{m}^{\lambda+\mu})\big)\Big)&
\end{align*}
Set $d_N=R_N^{\lambda}y_NR_N^{\lambda}+\lambda \bar{I}-\lambda R_N^{\lambda}$. Then according to the definition of $R_N^{\lambda}$, we have $d_N \le \lambda\bar{I}$. It follows that \begin{align*}
y_N & =d_N+(y_N-\lambda\bar{I})(\bar{I}-R_N^{\lambda})+(\bar{I}-R_N^{\lambda})(y_N-\lambda\bar{I})R_N^{\lambda}\\
&  \leq
\lambda\bar{I}+(y_N-\lambda\bar{I})(\bar{I}-R_N^{\lambda})+(\bar{I}-R_N^{\lambda})(y_N-\lambda\bar{I})R_N^{\lambda}.
\end{align*}
 Therefore, by the property of distribution
function, we have that
\begin{align*}
\nu(R_{m-1}^{\lambda+\mu}-R_{m}^{\lambda+\mu})
&\le\nu\Big({\bar{I}}_{[\mu,\infty)}\big((R_{m-1}^{\lambda+\mu}-R_{m}^{\lambda+\mu})\bar{\mathcal{E}}_m((y_N-\lambda\bar{I})
(\bar{I}-R_N^{\lambda})\\
&\quad\quad\quad+(\bar{I}-R_N^{\lambda})(y_N-\lambda\bar{I})R_N^{\lambda})(R_{m-1}^{\lambda+\mu}-R_{m}^{\lambda+\mu})\big)\Big).
\end{align*}
Using Chebyshev's inequality,  the right-hand side of the above expression is dominated by
\begin{align*}
&\frac{1}{{\mu}^2} \nu\Big(\big((R_{m-1}^{\lambda+\mu}-R_{m}^{\lambda+\mu})((y_N-\lambda\bar{I})(\bar{I}-R_N^{\lambda})\\
&\quad\quad\quad+(\bar{I}-R_N^{\lambda})(y_N-\lambda\bar{I})R_N^{\lambda})(R_{m-1}^{\lambda+\mu}-R_{m}^{\lambda+\mu})\big)^2\Big) \\
&\le \frac{1}{{\mu}^2} \nu\Big((R_{m-1}^{\lambda+\mu}-R_{m}^{\lambda+\mu})\big((y_N-\lambda\bar{I})(\bar{I}-R_N^{\lambda})\\
&\quad\quad\quad+(\bar{I}-R_N^{\lambda})(y_N-\lambda\bar{I})R_N^{\lambda}\big)^2(R_{m-1}^{\lambda+\mu}-R_{m}^{\lambda+\mu})\Big)\\
&=\frac{1}{{\mu}^2} \nu\big((R_{m-1}^{\lambda+\mu}-R_{m}^{\lambda+\mu})(y_N-\lambda \bar{I})
(\bar{I}-R_N^{\lambda})(y_N-\lambda \bar{I})(R_{m-1}^{\lambda+\mu}-R_{m}^{\lambda+\mu})\big)&\\
&+\frac{1}{{\mu}^2} \nu\big((R_{m-1}^{\lambda+\mu}-R_{m}^{\lambda+\mu})(\bar{I}-R_N^{\lambda})(y_N-\lambda\bar{I}) R_N^{\lambda}(y_N-\lambda\bar{I})(\bar{I}-R_N^{\lambda})(R_{m-1}^{\lambda+\mu}-R_{m}^{\lambda+\mu})\big)
\end{align*}
Summing $m$ from $n$ to $N$, we have that
\begin{eqnarray*}
&&\quad \nu (\bar{I} -R_N^{\lambda +\mu}) = \sum_{m=n}^N \nu (R_{m-1}^{\lambda+\mu}-R_{m}^{\lambda+\mu}) \\
&& \leq \frac{1}{{\mu}^2} \nu \big((\bar{I}-R_N^{\lambda+\mu })(y_N - \lambda \bar{I})(\bar{I} -R_N^{\lambda})(y_N - \lambda \bar{I})\big) \\
&& \ \ \  + \frac{1}{{\mu}^2} \nu \big((\bar{I}-R_N^{\lambda+\mu })(\bar{I}-R_N^{\lambda})(y_N - \lambda \bar{I})^2(\bar{I} -R_N^{\lambda})\big) \\
&&  \le  \frac{2}{{\mu}^2} \nu \big((\bar{I}-R_N^{\lambda}){(y_N- \lambda \bar{I})}^2\big). \\
\end{eqnarray*}
This completes the proof.
\end{proof}

The second one is an observation which will be quite essential in obtaining $e_{n+1,n+1}+e_{n+3,n+3}$ on the right hand side of \eqref{pro1}. This will in turn enable us to
deduce the desired distribution function inequality {{} for $x$ from \eqref{pro1} for $y$.}

\begin{lemma} \label{lem3.1} Let $\lambda>0$ and $z$ be a  positive operator in  $L_1(\mathcal{M})$.  Then for any  $m, k\geq n$, we have
$$  \nu (R_m^{\lambda}(e_{k+3,k+3} \otimes 1 \otimes z))=0.$$
\end{lemma}
\begin{proof}
By the definition of $y_n$, we have the following fact
\begin{align*}
&\quad y_n (e_{k+3,k+3}\otimes 1 \otimes z)\\
&=\big((e_{n+1,n+2} +e_{n+2,n+1}) \otimes \varepsilon_n \otimes \frac{(P|dx_{n}|^2 P)^\frac{1}{2}} {{\tau(P)}^\beta}\big) (e_{k+3,k+3} \otimes 1 \otimes z)\\
&=0 \otimes \varepsilon_n \otimes \big(\frac{(P|dx_{n}|^2 P)^\frac{1}{2}} {{\tau(P)}^\beta} z\big)\\
&=0
\end{align*}
which implies that  $R_n^{\lambda}(e_{k+3,k+3} \otimes 1 \otimes z)=0.$ Thus by the inequality $R_m^{\lambda} \le R_n^{\lambda}$ and $z\geq0$, we get that
$$
\nu (R_m^{\lambda}(e_{k+3,k+3} \otimes 1 \otimes z))\le \nu( R_n^{\lambda}(e_{k+3,k+3} \otimes 1 \otimes z))=0.
$$ This completes the proof.
\end{proof}

Now let us prove Proposition \ref{lem3.3}.

\begin{proof}
By Lemma \ref{lem3.2}, we just need to show $$ \nu \big((\bar{I}-R_N^{\lambda}){(y_N- \lambda \bar{I})}^2\big)\leq 2\left \| x \right \| ^2_{\mathcal{L}^c_{\beta}} \nu \big((\bar{I}-R_N^{\lambda})((e_{n+1,n+1}+e_{n+3,n+3}) \otimes 1 \otimes P)\big).$$

\textbf{Step 1.} Observe that
\begin{align*}
& \nu \big((\bar{I}-R_N^{\lambda}){(y_N- \lambda \bar{I})}^2\big) =  \sum_{m=n}^{N}
\nu\big((R_{m-1}^{\lambda}-R_{m}^{\lambda})(y_{N}-\lambda \bar{I} )R_{m}^{\lambda} (y_{N}- \lambda  \bar{I})\big)  \\
&\ \ \ \ \ \ \ \ \ \ \ \ \ \ \ \ \ \ \ \ \ \  \ \ \ \ \ \ \ \ \ \ \  \ \ +  \sum_{m=n}^{N}
\nu\big((R_{m-1}^{\lambda}-R_{m}^{\lambda})(y_{N}-\lambda \bar{I} )(\bar{I} -R_{m}^{\lambda})(y_{N}- \lambda  \bar{I})\big).
\end{align*}

Now write $\bar{I} -R_{m}^{\lambda}=\sum_{k=n}^{m}(R_{k-1}^{\lambda} -R_{k}^{\lambda})$ and  the last line of  the above expression is equal to
$$  \sum_{m=n}^{N} \sum_{k=n}^{m}
\nu\big((R_{m-1}^{\lambda}-R_{m}^{\lambda})(y_{N}-\lambda \bar{I} )(R_{k-1}^{\lambda} -R_{k}^{\lambda})(y_{N}- \lambda  \bar{I})\big).$$
Thus changing the order of
summation, we get
\begin{align*}
& \quad \nu \big((\bar{I}-R_N^{\lambda}){(y_N- \lambda \bar{I})}^2\big)\\
& =   \sum_{m=n}^{N}
\nu\big((R_{m-1}^{\lambda}-R_{m}^{\lambda})(y_{N}-\lambda \bar{I} )R_{m}^{\lambda}(y_{N}- \lambda  \bar{I})\big) \\
&\ \ \  +  \sum_{k=n}^{N}
\nu\big((R_{k-1}^{\lambda}-R_{N}^{\lambda})(y_{N}-\lambda \bar{I} )(R_{k-1}^{\lambda} -R_{k}^{\lambda})(y_{N}- \lambda  \bar{I})\big)  \\
& \le \sum_{m=n}^{N}
\nu\big((R_{m-1}^{\lambda}-R_{m}^{\lambda})(y_{N}-\lambda \bar{I} )(R_{m}^{\lambda} +R_{m-1}^{\lambda})(y_{N}- \lambda  \bar{I})\big).
\end{align*}
Using the martingale property of $ y $, we can split the above expression into two parts:
\begin{align*}
 &\quad\sum_{m=n}^{N}
\nu\big((R_{m-1}^{\lambda}-R_{m}^{\lambda})(y_{N}-y_m )(R_{m}^{\lambda} +R_{m-1}^{\lambda})(y_{N}-y_m )\big) \\
&\ \ \  +   \sum_{m=n}^{N}
\nu\big((R_{m-1}^{\lambda}-R_{m}^{\lambda})(y_m-\lambda \bar{I} )(R_{m}^{\lambda} +R_{m-1}^{\lambda})( y_m - \lambda  \bar{I})\big).
\end{align*}
Therefore, using $ (R_{m-1}^{\lambda} -R_m^{\lambda})y_m R_m^{\lambda} =0$ by the commuting property of $R^\lambda$, we get that
\begin{align*}
\begin{split}
\quad \nu \big((\bar{I}-R_N^{\lambda}){(y_N- \lambda \bar{I})}^2\big)
& \le 2 \sum_{m=n}^{N} \sum_{k=m+1}^{N}
\nu\big((R_{m-1}^{\lambda}-R_{m}^{\lambda})dy_k ^2\big)\\
& \ \ \ + \sum_{m=n}^{N}
\nu\big((R_{m-1}^{\lambda}-R_{m}^{\lambda})(y_m-\lambda \bar{I} )(R_{m-1}^{\lambda}-R_{m}^{\lambda} )( y_m - \lambda  \bar{I})\big) \\
& \le   2\sum_{m=n}^{N} \sum_{k=m+1}^{N}
\nu\big((R_{m-1}^{\lambda}-R_{m}^{\lambda})dy_k^2\big)\\
& \ \ \ +  \sum_{m=n}^{N}
\nu\big((R_{m-1}^{\lambda}-R_{m}^{\lambda})dy_m(R_{m-1}^{\lambda} -R_{m}^{\lambda})dy_m\big).\\
\end{split}
 \end{align*}
The last inequality follows from the following equivalent estimate
 $$ \nu\big((R_{m-1}^{\lambda}-R_{m}^{\lambda})(\lambda\bar{I}-y_{m-1})(R_{m-1}^{\lambda}-R_m^{\lambda})(2y_m-y_{m-1}-\lambda\bar{I})\big) \ge 0,$$
which in turn is deduced as follows: {by the definition of $R^\lambda$, $$ \nu\big((R_{m-1}^{\lambda}-R_{m}^{\lambda})(\lambda\bar{I}-y_{m-1})(R_{m-1}^{\lambda}-R_m^{\lambda})\big) \ge 0$$
 and \begin{align*}
  & \quad \nu\big((R_{m-1}^{\lambda}-R_{m}^{\lambda})(2y_m-y_{m-1}-\lambda\bar{I})\big)\\
  &  =2\nu\big((R_{m-1}^{\lambda}-R_{m}^{\lambda})(y_m-\lambda\bar{I})\big)+\nu\big((R_{m-1}^{\lambda}-R_{m}^{\lambda})(\lambda\bar{I}-y_{m-1})\big)\geq0.
  \end{align*}}
Hence we obtain  that
$$ \nu \big((\bar{I}-R_N^{\lambda}){(y_N- \lambda \bar{I})}^2\big) \le 2 \sum_{m=n}^{N} \sum_{k=m+1}^{N}
\nu\big((R_{m-1}^{\lambda}-R_{m}^{\lambda})dy_k^2\big)+ \sum_{m=n}^{N}
\nu\big((R_{m-1}^{\lambda}-R_{m}^{\lambda})dy_m^2\big).$$

\textbf{Step 2.} Note that \begin{align*}
&\quad  2\sum_{m=n}^{N} \sum_{k=m+1}^{N}
\nu\big((R_{m-1}^{\lambda}-R_{m}^{\lambda})dy_k^2\big)+ \sum_{m=n}^{N}
\nu\big((R_{m-1}^{\lambda}-R_{m}^{\lambda})dy_m^2\big)\notag\\
 & =2\sum_{m=n+1}^N \sum_{k=m+1}^{N} \nu \Big((R_{m-1}^\lambda-R_m^{\lambda})\big((e_{n+1,n+1}+e_{k+2,k+2})\otimes 1 \otimes \frac{P|dx_{k}|^2 P} {{\tau(P)}^{2\beta}}\big)\Big)\notag\\
 &\ \ \ +2 \sum_{k=n+1}^{N}\nu \Big((\bar{I}-R_n^{\lambda})\big((e_{n+1,n+1}+e_{k+2,k+2})\otimes 1 \otimes \frac{P|dx_{k}|^2 P} {{\tau(P)}^{2\beta}}\big)\Big)\notag\\
 &\ \ \ +\sum_{m=n+1}^{N} \nu \Big((R_{m-1}^{\lambda}-R_{m}^{\lambda})\big( (e_{n+1,n+1}+e_{m+2,m+2}) \otimes 1 \otimes \frac{P|dx_{m}|^2 P} {{\tau(P)}^{2\beta}}\big)\Big)\notag\\
 &\ \ \ +\nu \Big((\bar{I}-R_n^{\lambda}) \big((e_{n+1,n+1}+e_{n+2,n+2})  \otimes 1 \otimes \frac{P|dx_{n}|^2 P} {{\tau(P)}^{2\beta}} \big)\Big).\notag
\end{align*}
Then by Lemma \ref {lem3.1}, the above expression is equal to
\begin{align*}
&\quad 2\sum_{m=n+1}^N \sum_{k=m+1}^{N} \nu \big((R_{m-1}^\lambda-R_m^{\lambda})(e_{n+1,n+1}\otimes 1 \otimes \frac{P|dx_{k}|^2 P} {{\tau(P)}^{2\beta}})\big)\notag\\
 &\ \ \ +2 \sum_{k=n+1}^{N}\nu \big((\bar{I}-R_n^{\lambda})((e_{n+1,n+1}+e_{k+2,k+2})\otimes 1 \otimes \frac{P|dx_{k}|^2 P} {{\tau(P)}^{2\beta}})\big)\notag\\
 &\ \ \ +\sum_{m=n+1}^{N} \nu \big((R_{m-1}^{\lambda}-R_{m}^{\lambda})( e_{n+1,n+1} \otimes 1 \otimes \frac{P|dx_{m}|^2 P} {{\tau(P)}^{2\beta}})\big)\notag\\
 &\ \ \ +\nu \big((\bar{I}-R_n^{\lambda}) ((e_{n+1,n+1}+e_{n+2,n+2})  \otimes 1 \otimes \frac{P|dx_{n}|^2 P} {{\tau(P)}^{2\beta}} )\big).\notag
\end{align*}
Now, grouping together the terms involving $e_{n+1,n+1}$, we get that
\begin{eqnarray}\label{dis fff}
&& \quad  2\sum_{m=n}^{N} \sum_{k=m+1}^{N}
\nu\big((R_{m-1}^{\lambda}-R_{m}^{\lambda})dy_k^2\big)+ \sum_{m=n}^{N}
\nu\big((R_{m-1}^{\lambda}-R_{m}^{\lambda})dy_m^2\big)\notag\\
&& \leq2\sum_{m=n}^N \nu \big((R_{m-1}^\lambda-R_m^{\lambda})(e_{n+1,n+1}\otimes 1 \otimes\frac{P( S_c^2 (x_{N}-x_{m-1}) P} {{\tau(P)}^{2\beta}})\big)\notag\\
 &&\ \ \ +2 \sum_{k=n+1}^{N}\nu ( e_{k+2,k+2}  \otimes 1 \otimes \frac{P|dx_{k}|^2 P} {{\tau(P)}^{2\beta}} )\notag\\
 &&\ \ \ +\nu \big((\bar{I}-R_n^{\lambda}) (e_{n+2,n+2}  \otimes 1 \otimes \frac{P|dx_{n}|^2 P} {{\tau(P)}^{2\beta}} )\big).
\end{eqnarray}
Observe that by Lemma  \ref {lem3.1}, the trace-preserving property of conditional expectations and $R_N^{\lambda}\leq R_n^{\lambda}$,
\begin{eqnarray}\label{dis ggg}
&&\quad 2 \sum_{k=n+1}^{N}\nu ( e_{k+2,k+2}  \otimes 1 \otimes \frac{P|dx_{k}|^2 P} {{\tau(P)}^{2\beta}} )\notag\\
 &&+\nu \big((\bar{I}-R_n^{\lambda}) (e_{n+2,n+2}  \otimes 1 \otimes \frac{P|dx_{n}|^2 P} {{\tau(P)}^{2\beta}} )\big)\notag\\
&&\leq 2\sum_{k=n}^{N}\tau ( \frac{P|dx_{k}|^2 P} {{\tau(P)}^{2\beta}}) \notag\\
&&= 2\nu \big((\bar{I}-R_n^{\lambda}){\bar{\mathcal{E}}}_{n}(e_{n+3,n+3} \otimes 1 \otimes  \frac{PS_c^2(x_{N}-x_{n-1}) P} {{\tau(P)}^{2\beta}})\big)\notag\\
&& \leq 2\nu \big((\bar{I}-R_N^{\lambda})(e_{n+3,n+3} \otimes 1 \otimes  \frac{P\mathcal{E}_n(S_c^2(x_{N}-x_{n-1})) P} {{\tau(P)}^{2\beta}})\big).
\end{eqnarray}
Putting (\ref{dis fff}) and (\ref{dis ggg}) together, by the trace-preserving property again and $ R_N^{\lambda}\leq  R_n^{\lambda}$, we see that
 \begin{eqnarray}\label{dis aaa}
&&\quad  2\sum_{m=n}^{N} \sum_{k=m+1}^{N}
\nu\big((R_{m-1}^{\lambda}-R_{m}^{\lambda})dy_k^2\big)+ \sum_{m=n}^{N}
\nu\big((R_{m-1}^{\lambda}-R_{m}^{\lambda})dy_m^2\big)\notag\\
  && \leq 2\sum_{m=n}^N \nu \big((R_{m-1}^\lambda-R_m^{\lambda}) (e_{n+1,n+1}\otimes 1 \otimes\frac{P \mathcal{E}_m(S_c^2 (x_{N}-x_{m-1})) P} {{\tau(P)}^{2\beta}})\big)\notag\\
 &&\ \ \ +2\nu \big((\bar{I}-R_N^{\lambda})(e_{n+3,n+3} \otimes 1 \otimes  \frac{P{\mathcal{E}}_{n}(S_c^2(x_{N}-x_{n-1})) P} {{\tau(P)}^{2\beta}})\big)\notag\\
   && \leq 2\sum_{m=n}^N \nu \big((R_{m-1}^\lambda-R_m^{\lambda})(e_{n+1,n+1} \otimes 1
\otimes P)\big)\cdot \sup_{m\geq n} \frac{\| P {\mathcal{E}}_{m}(S_c^2 (x_{N}-x_{m-1})) P \|_{\infty}}{{{\tau(P)}^{2\beta}}}\notag \\
  && \ \ \  +2 \nu \big((\bar{I}-R_N^{\lambda})(e_{n+3,n+3} \otimes 1
\otimes P)\big)\cdot \frac{\| P {\mathcal{E}}_{m}(S_c^2 (x_{N}-x_{m-1})) P \|_{\infty}}{{{\tau(P)}^{2\beta}}}.
\end{eqnarray}
By duality,  we can write for any $k\geq n$
 \begin{eqnarray*}\label{dis aab}
 &&\frac{{\| P {\mathcal{E}}_{k}(S_c^2 (x_{N}-x_{k-1})) P \|_{\infty}}}{{{\tau(P)}^{2\beta}}} \notag\\
 && =\sup_{\|z\|_1\leq1, \ z\in L_1^+(P\mathcal{M}_{k}P)}\frac{\tau(zPS_c^2(x_{N}-x_{k-1})P)}{\tau(P)^{2\beta}}\notag\\
 && =\sup_{P'\in{\mathcal{P}(P\mathcal{M}_{k}P)}}\frac{1}{\tau(P')}\frac{\tau(P'PS_c^2(x_{N}-x_{k-1})P)}{\tau(P)^{2\beta}},\notag\\
 \end{eqnarray*}
 where in the last equality we have exploited the extreme point property of $L_1$, that is, any $z\in L_1(P\mathcal{M}_{k}P)$ can be rewritten as $$z=\sum_k\lambda_k\frac{e_k}{\tau(e_k)}$$ with $e_k$'s projections and
 $\|z\|_1=\sum_k|\lambda_k|$.
It follows that
\begin{align}\label{dis aac}
\frac{\|P {\mathcal{E}}_{k}(S_c^2 (x_{N}-x_{k-1})) P \|_\infty} {{\tau(P)}^{2\beta}}
\leq {\sup_{k\geq1}\sup_{P'\in{\mathcal{P}(P\mathcal{M}_{k}P)}}{\frac{\|(x-x_{k-1})P'\|_2^2}{\tau(P')^{2\beta+1}}}}\leq \|x\|_{\mathcal{L}^c_{\beta}}^2
\end{align}
 since  $P'\leq P$ and
 $$\tau(P'S_c^2(x_{N}-x_{k-1})P')=\|(x_N-x_{k-1})P'\|_2^2.$$

 Putting (\ref{dis aaa}) and (\ref{dis aac}) together, we get that
\begin{eqnarray*}\label{dis aad}
&&\quad 2 \sum_{m=n}^{N} \sum_{k=m+1}^{N}
\nu\big((R_{m-1}^{\lambda}-R_{m}^{\lambda})dy_k^2\big)+ \sum_{m=n}^{N}
\nu\big((R_{m-1}^{\lambda}-R_{m}^{\lambda})dy_m^2\big)\notag\\
&&\leq 2\|x\|_{\mathcal{L}^c_{\beta}}^2 \nu \big((\bar{I}-R_N^{\lambda})((e_{n+1,n+1}+e_{n+3,n+3}) \otimes 1
\otimes P)\big).\notag\\
\end{eqnarray*}
Thus we obtain the desired estimate.
\end{proof}

With Proposition \ref{lem3.3}, we are now at the position to conclude John-Nirenberg theorem in terms of distribution function inequality.

\begin{theorem}\label{lem3.4}   Let $\beta\geq0$ and $x\in \mathcal L^c_\beta(\mathcal M)$. Then for any $n\geq1$, $P\in\mathcal P(\mathcal M_n)$ and $\lambda\geq0$, there holds
\begin{align}\label{lem3.40} \tau \left({{I}}_{[\lambda,\infty)}\Big(\frac{(P S_c^2 (x-x_{n-1}) P)^{\frac{1}{2}}} {{\tau(P)}^{\beta}}\Big)\right)  \le 2(1-e^{-2})^{-1} e^2 e^{- \frac{ \lambda}{e{\| x \| }_{\mathcal{L}_{\beta}^c}}} \tau(P).\end{align}
\end{theorem}
\begin{proof}
Let $x\in \mathcal{L}_{\beta}^c(\mathcal M)$. Fixed one $n\geq1$ and $N\geq n$.
For a fixed positive integer $k$, apply Proposition \ref{lem3.3} inductively with $\lambda = (k-1)e,\dotsm,e,0$ and  $\mu=e$,
\begin{eqnarray*}
\nu(\bar{I}-R_{N}^{ke}) & \le &\frac{4}{e^2} \|x\|^2_{\mathcal{L}^c_{\beta}} \nu \big((\bar{I}-R_N^{(k-1)e})((e_{n+1,n+1}+e_{n+3,n+3} )\otimes 1 \otimes P)\big)\\
 & \le &\frac{4}{e^2} \|x \|^2_{\mathcal{L}^c_{\beta}} \nu (\bar{I}-R_N^{(k-1)e})\\
& \le & (\frac{4}{e^2} \|x \|^2_{\mathcal{L}^c_{\beta}})^{k} \nu \big((\bar{I}-R_N^0)((e_{n+1,n+1}+e_{n+3,n+3} )\otimes 1 \otimes P)\big) \\
 & \le &(\frac{4}{e^2} \|x \|^2_{\mathcal{L}^c_{\beta}})^{k} \nu \big((e_{n+1,n+1}+e_{n+3,n+3} )\otimes 1 \otimes P\big) \\
  & = &2(\frac{4}{e^2} \|x \|^2_{\mathcal{L}^c_{\beta}})^{k} \tau (P). \\
\end{eqnarray*}
Let $m \ge n ,$ set $P_m^{ke}= \bigwedge\limits_{l \ge k} R_m^{le}$ for any  $ k \in \mathbb{N},$ then one observes easily that  $P_m^{le} \le P_m^{ke} $ if $l \le k.$  By the definition of $P_N^{ke}$, one may relate $P_N^{ke}$ with $R_N^{(k+1)e}$ as follows,
\begin{eqnarray*}
\nu(\bar{I}-R_N^{ke}) & = &\nu(\bar{I}-P_N^{ke}) +\nu(P_N^{ke}-R_N^{ke})\\
& = & \nu(\bar{I}-P_N^{ke}) +\nu(R_N^{ke} \wedge P_N^{(k+1)e} - R_N^{ke})\\
& = & \nu(\bar{I}-P_N^{ke}) -\nu(R_N^{ke} - R_N^{ke} \wedge P_N^{(k+1)e} )\\
& = & \nu(\bar{I}-P_N^{ke}) -\nu(  R_N^{ke} \vee P_N^{(k+1)e} -P_N^{(k+1)e} )\\
& \ge & \nu(\bar{I}-P_N^{ke}) -\nu(  \bar{I}- P_N^{(k+1)e} )=\nu(  P_N^{(k+1)e} -P_N^{ke} ),
\end{eqnarray*}
where in last equality, we have used the fact that for any two projections $P, Q$, there holds the equivalence
\begin{align}\label{equiv}P-P\wedge Q\sim P\vee Q-Q.\end{align}

Assume that $ \| x \|_{\mathcal{L}^c_{\beta}} = \frac{1}{2}$.  Putting the  two inequalities above together, we get that
\begin{eqnarray}\label{dis aae}
\nu(  P_N^{(k+1)e}-P_N^{ke} ) \le \nu(  \bar{I}-R_N^{ke} )  \le 2 e^{-2k} \tau(P).
\end{eqnarray}

Now fix one $ \lambda \geq0$. Then there exists  a unique positive integer $k_0$ such that  $k_0e \le \lambda < (k_0+1)e$.
We claim first that
$ I_{[k_0e,\infty)}(y_N)$ is equivalent to
a subprojection of $\bar{I}-P_N^{k_0e}$. Indeed,
for any nonzero vector  $\xi \in P_N^{k_0e} (H),$ since $P_N^{k_0e}\leq R_N^{k_0e}$ and $R_N^{k_0e}y_NR_N^{k_0e}< k_0e$ by definitions, we obtain
$$ (y_N \xi ,\xi) =(P_N^{k_0e} y_N P_N^{k_0e} \xi, \xi) < k_0e \left \| \xi  \right \| ^2$$ which implies that
 $\xi \notin I_{[k_0e,\infty)}(y_N)(H)$. In other words,
 $$P_N^{k_0e}\wedge  \bar{I}_{[k_0e,\infty)}(y_N)=0.$$ Then by \eqref{equiv}, one concludes the claim,
 \begin{align*}
  \bar{I}_{[k_0e,\infty)}(y_N)&= \bar{I}_{[k_0e,\infty)}(y_N)- \bar{I}_{[k_0e,\infty)}(y_N)\wedge P_N^{k_0e}\\
  &\sim\bar{I}_{[k_0e,\infty)}(y_N)\vee P_N^{k_0e}-P_N^{k_0e}\\
  &\leq \bar{I}-P_N^{k_0e}.
 \end{align*}
 Therefore, by (\ref{dis aae}), we deduce that
\begin{eqnarray*}
\nu({\bar{I}}_{[\lambda,\infty)}( y_N)) & \le & \nu(\bar{I}_{[k_0e,\infty)}( y_N))\\
&\le & \nu (\bar{I}-P_N^{k_0e})=\begin{matrix}
\sum_{l=k_0}^{\infty} \nu(  P_N^{(l+1)e} -P_N^{le} )
\end{matrix} \\
& \le & 2\tau(P)\begin{matrix}
\sum_{l=k_0}^{\infty} e^{-2l}
\end{matrix} =2\tau(P)(1-e^{-2})^{-1} e^{-2k_0}.\\
& \leq & 2 \tau(P)(1-e^{-2})^{-1} e^{(2-\frac{2\lambda}{e})}.
\end{eqnarray*}
By homogeneity, we have for $x\in \mathcal L^c_\beta(\mathcal M)$ not necessarily of norm $1/2$,
\begin{align*}
\nu(\bar{I}_{[\lambda,\infty)}( y_N  )) \le 2(1-e^2)^{-1} e^2 e^{- \frac{ \lambda}{e{\| x \|}_{\mathcal{L}_{\beta}^c}}} \tau(P).
\end{align*}
Finally, noticing that $y_N^2 \ge  e_{n+1,n+1}\otimes 1 \otimes\frac{P S_c^2 (x_{N}-x_{n-1}) P}{{\tau(P)}^{2\beta}} $,  we arrive at
\begin{eqnarray*}
&&\quad \tau\big({{I}}_{[\lambda,\infty)}(\frac{(P S_c^2 (x_{N}-x_{n}) P)^{\frac{1}{2}}}{{\tau(P)}^{\beta}})\big)\\
&& =  \nu \big({\bar{I}}_{[\lambda,\infty)}(e_{n+1,n+1}\otimes 1 \otimes\frac{(P S_c^2 (x_{N}-x_{n}) P)^{\frac{1}{2}}}{{\tau(P)}^{\beta}} ) \big)\\
&&\le  \nu (\bar{I}_{[\lambda,\infty)}( y_N )) \le  2(1-e^{-2})^{-1} e^2 e^{- \frac{ \lambda}{e{\| x \|}_{\mathcal{L}_{\beta}^c}}} \tau(P).
\end{eqnarray*}
Let $N\rightarrow\infty$, this concludes the proof.
\end{proof}

\subsection{Exponential integrability and the $p$-moment characterization for $\mathcal {L}_\beta^c(\mathcal M)$}
\label{subsection 2.2}
With the help of Theorem \ref{lem3.4}, we obtain the  following  John-Nirenberg theorem in terms of exponential integrability.

\begin{theorem}\label{lem3.21} Let $\beta\geq0$ and $x\in \mathcal L^c_\beta(\mathcal M)$. Then for any $0<a<\frac{1 }{e{\left \| x \right \| }_{{\mathcal{L}}_{\beta}^c}}$,
$$\sup_{n\geq 1}\sup_{P \in {\mathcal{P}(\mathcal{M}_n)}}\frac1{\tau(P)}\tau\Big(e^{\frac{a(PS_c^2(x-x_{n-1})P)^{\frac{1}{2}}}{\tau(P)^{\beta}}}\Big)\leq K_a<\infty,$$
where $K_a=1+2a(1-e^{-2})^{-1} e^2\int_0^\infty e^{(a-\frac{1}{e{\|  x  \|}_{{\mathcal{L}}_{\beta}^c}})s}ds.$
\end{theorem}

\begin{proof}  It suffices to consider the case of $\|x\|_{{\mathcal{L}^c_{\beta}}}=1$ and  $0<a<\frac{1}{e}$. Fix $n\geq1$ and $P \in \mathcal{P}(\mathcal{M}_n)$.
  Set $\Phi(s)=e^{as}-1$. Then by the equality $\tau(\Phi(|y|))=\int_0^\infty \lambda_s(|y|)d\Phi(s)$ and Theorem \ref{lem3.4},
we have
\begin{eqnarray*}
\tau\Big(e^{\frac{a(PS_c^2(x-x_{n-1})P)^{\frac{1}{2}}}{\tau(P)^{\beta}}}-P\Big)
&=&\int_0^\infty \lambda_s\Big(\frac{(PS_c^2(x-x_{n-1})P)^{\frac{1}{2}}}{(\tau(P))^{\beta}}\Big)de^{as}\\
&\leq&a\int_0^\infty e^{as}\cdot\tau(P) 2(1-e^{-2})^{-1} e^2e^{\frac{-s}{e}}ds\\
&=&2a(1-e^{-2})^{-1} e^2\tau(P)\int_0^\infty e^{(a-\frac{1}{e})s}ds,
\end{eqnarray*}
which gives immediately the desired estimates.
\end{proof}

Next we consider the $p$-moment inequality form of John-Nirenberg theorem.

\begin{definition}
Let $\beta\geq0$ and $0< p<\infty$. We define
$$\mathcal{L}^c_{\beta,p}(\mathcal{M})=\{x\in L_2(\mathcal{M}): \|x\|_{\mathcal{L}^c_{\beta,p}}<\infty\}$$ where
$$
\|x\|_{\mathcal{L}^c_{\beta,p}}= \sup_{n\geq1}\sup_{e\in{\mathcal{P}(\mathcal{M}_n)}}{\frac{\|S_c(x-x_{n-1})e\|_{p}}{(\tau(e))^{\beta+\frac{1}{p}}}},
$$
and $$\mathcal{L}^r_{\beta,p}(\mathcal{M})=\{x\in L_2(\mathcal{M}): x^*\in \mathcal{L}^c_{\beta,p}(\mathcal{M})\}.$$
\end{definition}

\begin{remark}
(i) Note that when $p=2$ we have that $\mathcal{L}^c_{\beta,p}(\mathcal{M})=\mathcal{L}^c_{\beta}(\mathcal{M}).$

(ii) The above definition is motivated by the John-Nirenberg theorem for $\mathsf{bmo}^c(\mathcal M)$ obtained in \cite{ref5}. When $\beta=0$, it is not difficult to check that this definition goes back to the classical one given in \cite[Definition 2.45]{ref25}. For $\beta>0$, however, whenever $\mathcal M$ is commutative or not, this definition seems new in the literature.
\end{remark}

\begin{theorem}\label{lem3.25}  Let $\beta\geq0$, $x\in\mathcal{L}^c_{\beta}(\mathcal{M})$ and $0< p<\infty$.  Then there exist two constants
$\alpha_p$ and $\beta_p$ such that
\begin{eqnarray}\label{dis clx}
\alpha_p^{-1}\|x\|_{\mathcal{L}^c_{\beta}}\leq\|x\|_{\mathcal{L}^c_{\beta,p}}\leq\beta_p\|x\|_{\mathcal{L}^c_{\beta}}
\end{eqnarray}
with $\alpha_p$  and $\beta_p$  satisfying \ \
\begin{center} \ \ \ \ \ \  \  \rm{(i)} $\alpha_p=1$ for $2\leq p<\infty$, \ \ \   \rm{(ii)}  $\alpha_p\leq C^{1-2/p}$ for $0< p<2$,\end{center}
\begin{center} \rm{(iii)} $\beta_p\leq cp$ for $2\leq p<\infty$, \  \rm{(iv)} $\beta_p=1$ for $0< p<2$.\end{center}
\noindent{ The same inequalities hold for $\|\cdot\|_{\mathcal{L}^r_{\beta}}$ and $\|\cdot\|_{\mathcal{L}^r_{\beta,p}}$.}
\end{theorem}

\begin{proof}
We only need to prove the column case, since the row case can be done by replacing $x$
with $x^*$.  First we show the inequality on the right hand side.  For $2\leq p<\infty$, by Theorem \ref{lem3.4},
 we have that for any $n\geq1$ and $P\in \mathcal{P}(\mathcal{M}_n)$
\begin{eqnarray*}
{\frac{\|S_c(x-x_{n-1})P\|_{p}}{(\tau(P))^{\beta+\frac{1}{p}}}}
&=&(\tau(P))^{-\frac{1}{p}}\big(p\int_0^\infty{s}^{p-1}\cdot{\lambda_s\big(\frac{(PS_c^2(x-x_{n-1})P)^{\frac{1}{2}}}{(\tau(P))^{\beta}}\big)}ds\big)^{\frac{1}{p}}\notag\\
&\leq&\big(p\int_0^\infty{s}^{p-1}\cdot 2(1-e^{-2})^{-1} e^2 e^{- \frac{s}{e{ \| x \|}_{\mathcal{L}_{\beta}^c}}} ds\big)^{\frac{1}{p}}\notag\\
&=&e(2p(1-e^{-2})^{-1}e^2)^{\frac{1}{p}}\Gamma(p)^{\frac{1}{p}}\|x\|_{\mathcal{L}^c_{\beta}}.
\end{eqnarray*}
Noting  that $$\lim_{p\rightarrow\infty}\frac{(2p(1-e^2)^{-1}e^2)^{\frac{1}{p}}\Gamma(p)^{\frac{1}{p}}}{p}=\frac{1}{e},$$  and thus there exists a constant $c$ such that   for $2\leq p<\infty$ \begin{eqnarray}\label{dis cly}
\|x\|_{\mathcal{L}^c_{\beta,p}}\leq cp\|x\|_{\mathcal{L}^c_{\beta}}.\end{eqnarray}
The case $0<p<2.$  Let $q>0$ such that $1/p=1/q+1/2$. By H$\ddot{\mbox{o}}$lder's inequality, one obtains  for any $n\geq1$ and $P\in \mathcal{P}(\mathcal{M}_n)$
$$\|S_c(x-x_{n-1})P\|_{p}\leq(\tau(P))^{\frac{1}{q}}\|S_c(x-x_{n-1})P\|_{2}.$$
Thus we have
$$(\tau(P))^{-\beta-\frac{1}{p}}\|S_c(x-x_{n-1})P\|_{p}\leq(\tau(P))^{-\beta-\frac{1}{2}}\|S_c(x-x_{n-1})P\|_{2}$$
which implies that $$\|x\|_{\mathcal{L}^c_{\beta,p}}\leq\|x\|_{\mathcal{L}^c_{\beta,2}}=\|x\|_{\mathcal{L}^c_{\beta}}.$$

We turn to the inequality on the left hand side. First we consider the case of $2\leq p<\infty$.
Choose $q$ such that $1=2/p+1/q$. Fix $n$ and  $e\in \mathcal{P}(\mathcal{M}_n)$. By  H$\ddot{\mbox{o}}$lder's inequality,
one gets
 \begin{eqnarray*}
 \sup_{e\in{\mathcal{P}(\mathcal{M}_n)}}\frac{\|(x-x_{n-1})e\|_2^2}{(\tau(e))^{1+2\beta}}
&=&\sup_{e\in{\mathcal{P}(\mathcal{M}_n)}}\frac{{\tau(eS_c^2(x-x_{n-1})e)}}{(\tau(e))^{1+2\beta}}\\
&\leq&\sup_{e\in{\mathcal{P}(\mathcal{M}_n)}}\frac{1}{{(\tau(e))^{1+2\beta}}}(\|e\|_q\cdot\|e{S}_c^2(x-x_{n-1})e\|_{\frac{p}{2}})\\
&=&\sup_{e\in{\mathcal{P}(\mathcal{M}_n)}}{\frac{\|S_c(x-x_{n-1})e\|_{p}^2}{(\tau(e))^{2\beta+\frac{2}{p}}}}.\\
\end{eqnarray*}
Thus we have $\|x\|_{\mathcal{L}^c_{\beta}}\leq\|x\|_{\mathcal{L}^c_{\beta,p}}$.
 For $0< p<2$, fix one $q>2$,
choose $0<\theta<1$, such that $1=(2\theta)/p+(1-\theta)/q$.  By  H$\ddot{\mbox{o}}$lder's inequality,
one gets
 \begin{eqnarray*}
\frac{\|(x-x_{n-1})e\|_2^2}{(\tau(e))^{1+2\beta}}
&=&\frac{\tau\big((eS_c^2(x-x_{n-1})e)^{\theta+(1-\theta)}\big)}{(\tau(e))^{1+2\beta}}\\
&\leq& \frac{1}{(\tau(e))^{1+2\beta}}\|(eS_c^2(x-x_{n-1})e)^{\theta}\|_{\frac{p}{2\theta}}\|(eS_c^2(x-x_{n-1})e)^{1-\theta}\|_{\frac{q}{1-\theta}}\\
&=& \big(\frac{\|S_c(x-x_{n-1})e\|_{p}}
{(\tau(e))^{\beta+\frac{1}{p}}}\big)^{2\theta}\big(\frac{\|S_c(x-x_{n-1})e\|_{2q}}{(\tau(e))^{\beta+\frac{1}{2q}}}\big)^{2(1-\theta)}\\
\end{eqnarray*}
 which implies that\begin{eqnarray*}
\|x\|_{\mathcal{L}^c_{\beta}}\leq\|x\|^\theta_{\mathcal{L}^c_{\beta,p}}\|x\|^{1-\theta}_{\mathcal{L}^c_{\beta,2q}}.
\end{eqnarray*}
Thus by the upper inequalities in \eqref{dis clx}, we have
 \begin{eqnarray*}
\|x\|_{\mathcal{L}^c_{\beta}}\leq (cq)^{\frac{1-\theta}{\theta}}\|x\|_{\mathcal{L}^c_{\beta,p}}.
\end{eqnarray*}
Noting that ${\frac{1-\theta}{\theta}}=\frac{1-\frac{2}{p}}{\frac{1}{q}-1}$, we get the desired estimate by taking
$C=(cq)^{1/(1/q-1)}$.
\end{proof}

\begin{remark}\label{lem3.8}
(i) Theorem \ref{lem3.25} (or its proof) actually tells us that the Lipschitz space $\mathcal L^c_{\beta}(\mathcal M)$ coincides with $\mathcal L^c_{\beta,p}(\mathcal M)$ for any $2\leq p<\infty$. While for $0< p<2$, if a priori we assume that $x\in\mathcal L^c_{\beta}(\mathcal M)$, then the norms are equivalent; so it would be interesting to show a distribution inequality as \eqref{lem3.40}  starting with $x\in \mathcal L^c_{\beta,p}(\mathcal M)$.

(ii) Another closely related question is whether there exists a direct approach to the John-Nirenberg theorem for mixed Lipschitz spaces via the distribution function inequality; an indirect way via duality and interpolation has been provided in \cite[Theorem 3.20]{ref5} for mixed BMO space.

\end{remark}

\subsection{John-Nirenberg theorem for $\Lambda^c_\beta(\mathcal M)$}
\label{subsection 2.3}

The above proof for $\mathcal L^c_\beta(\mathcal M)$ works equally for $\Lambda^c_\beta(\mathcal M)$. Indeed, for $x\in \Lambda^c_\beta(\mathcal M)$, it suffices to modify the definition  in (\ref{sss}) as $$ y_{m} = \sum_{k=n}^m(e_{n+1,k+2} + e_{k+2,n+1 })  \otimes \varepsilon_{k} \otimes \frac{(P{\mathcal{E}_{k} }|dx_{k+1}|^2 P)^\frac{1}{2}} {{\tau(P)}^\beta}, $$
and follow the rest of the arguments verbatim, one may obtain the distribution function inequality and thus the exponential inequality as well as the $p$-moment characterization. We {leave} the details to the interested reader, and just summarize the results in the following theorem. For the Hardy spaces $\mathsf{h}^c_p(\mathcal M)$ appearing below, we refer the reader to next section for more information.

\begin{theorem}\label{thm:JNlambda}   Let $\beta\geq0$ and $x\in \Lambda^c_\beta(\mathcal M)$. The following assertions hold.

{\rm (i)} For any $n\geq1$, $P\in\mathcal P(\mathcal M_n)$ and $\lambda\geq0$, there holds
\begin{align}\label{JNlambdaDistribution} \tau \left({{I}}_{[\lambda,\infty)}\Big(\frac{(Ps_c^2 (x-x_{n}) P)^{\frac{1}{2}}} {{\tau(P)}^{\beta}}\Big)\right)  \le 2(1-e^{-2})^{-1} e^2 e^{- \frac{ \lambda}{e{\| x \| }_{{\Lambda_{\beta}^c}}}} \tau(P).
\end{align}

{\rm (ii)} For any $0<a<\frac{1 }{e{\| x\| }_{\Lambda_{\beta}^c}}$,
$$\sup_{n\geq 1}\sup_{P \in {\mathcal{P}(\mathcal{M}_n)}}\frac1{\tau(P)}\tau\Big(e^{\frac{a(Ps_c^2(x-x_{n})P)^{\frac{1}{2}}}{\tau(P)^{\beta}}}\Big)\leq K_a,$$
{{}where $K_a$ is the constant appeared in Theorem \ref{lem3.21}.}

{\rm (iii)}
Let $0< p<\infty$.  There holds
\begin{eqnarray}\label{JNlambdaMoment}
\alpha_p^{-1}\|x\|_{\Lambda^c_{\beta}}\leq\|x\|_{\Lambda^c_{\beta,p}}\leq\beta_p\|x\|_{\Lambda^c_{\beta}},
\end{eqnarray}
where $\alpha_p,\beta_p$ are the constants appeared in Theorem \ref{lem3.25} and
$$
\|x\|_{{\Lambda}^c_{\beta,p}}=\max\{\|\mathcal{E}_1(x)\|_\infty,\ \ \sup_{n\geq1}\sup_{e\in{\mathcal{P}(\mathcal{M}_n)}}{\frac{\|(x-x_n)e\|_{{\mathsf{h}}_p^c}}{(\tau(e))^{\beta+\frac{1}{p}}}}\}.
$$
\end{theorem}


\begin{remark}\label{lem3.5}
(i) Note that the characterization \eqref{JNlambdaMoment} for $1\leq p<\infty$ without explicit estimates over $\alpha_p$ and $\beta_p$ has been obtained in \cite{ref3} as a corollary of their atomic decomposition and duality results. While our direct method gives the order $\beta_p\leq c{p}$ for $2\leq p<\infty$, and extends the scale of $p$ {{}in \eqref{JNlambdaMoment} to $0<p<1$.}

(ii) When $\beta=0$, we recover the $p$-moment form of John-Nirenberg theorem for ${\mathsf{bmo}}^c$ spaces (see \cite{ref5}).
\end{remark}

\section{Symmetric space moment characterization for  $\Lambda^c_\beta(\mathcal M)$}
\label{section 3}
This section is devoted to the symmetric space moment characterization of $\Lambda^c_\beta(\mathcal M)$. For this purpose, we need some properties of symmetric Banach function spaces and noncommutative symmetric spaces.

\subsection{Symmetric Banach function spaces on $(0,\infty)$}
\label{subsection 3.1}
Denote by  $L_0(0,\infty)$ be the space of all Lebesgue measurable
real-valued functions  on $(0,\infty) $. Let $f\in L_0(0,\infty)$. Recall that the distribution function of $f$ is defined by
$$\lambda_s(f)=|\{\omega\in (0,\infty): |f(\omega)|>s\}|, \ s>0$$
and its non-increasing rearrangement by $$\mu_t(f)=\inf\{s>0:\lambda_s(f)\leq t\}, \ \ t>0.$$
 A quasi-Banach subspace $E$ of $L_0(0,\infty)$ is called symmetric if for any
$g\in E$ and any measurable function $f$ with $\mu_t(f)\leq \mu_t(g)$ for all $t\geq0$, then $f\in E$ and $\|f\|_E\leq \|g\|_E$.
A symmetric quasi-Banach space $E$ is said to have the Fatou property if for every net $(f_i)_{i\in I}$ in $E$
satisfying $0\leq f_i \uparrow$ and $\sup_{i\in I} \|f_i\|_E<\infty$ the supremum $f=\sup_{i\in I}f_i$ exists in $E$ and $\|f_i\|_E\uparrow
\|f\|_E$. We say that $E$ has order continuous norm if for every net $(f_i)_{i\in I}$ in $E$ such that $f_i \downarrow0$ we have  $\|f_i\| \downarrow0$.

The
K${\ddot{\mbox{o}}}$the dual of a symmetric Banach space $E$  is  given by
$$E^{\times}=\{f\in L_0(0,\infty) :\int^\infty_0|f(t)g(t)|dt<\infty,\; \forall g\in E \},$$
with the norm
$\|f\|_{E^\times}:=\sup\{\int^\infty_0|f(t)g(t)|dt:\|g\|_E\leq1\}.$ The space $E^{\times}$ is  symmetric and has the Fatou property.
A symmetric Banach function space $E$ on $(0,\infty)$ has order continuous norm if and only if it is separable, which is also equivalent to the statement $E^*=E^\times$.  We refer to \cite{ref1,ref19,ref4} for more details.

 For any $s>0$,  we define the dilation operator $D_s$ on $L_0(0,\infty)$ by
$$(D_sf)(t)=
f(t/s),\ \  t>0.$$
 For a quasi-Banach function space $E$, the
lower and upper Boyd indices $p_E$ and $q_E$ of $E$  are respectively defined by
$$p_E:=\lim\limits_{s\rightarrow\infty}\frac{\log s}{\log\|D_s\|}  \ \ \ \ \   \mbox{and} \ \ \ \ \
q_E:=\lim\limits_{s\rightarrow 0^{+}}\frac{\log s}{\log\|D_s\|}.$$
$D_s$ is a bounded linear operator on $E$ for every $s>0$ and $0\leq p_E\leq q_E\leq\infty$.
If $E$ is furthermore a Banach space, then  $1\leq p_E\leq q_E\leq\infty$ and
\begin{align}\label{dis cll}\frac{1}{p_E}+\frac{1}{q_{E^{\times}}}=1,\ \  \frac{1}{p_{E^{\times}}}+\frac{1}{q_{E}}=1.\end{align}
Note that if $E$ is a separable symmetric Banach space with  $1< p_E\leq q_E<\infty$, then $E$ automatically have the Fatou property.

Given a quasi-Banach function space $E$, for $0<r<\infty$, $E^{(r)}$ will denote the quasi-Banach function space
defined by $E^{(r)}=\{f: |f|^r\in E\}$
equipped with the quasi-norm $\|f\|_{E^{(r)}}=\big\||f|^r\big\|^{\frac{1}{r}}_E.$
Note that\begin{align}\label{dis clr}p_{E^{(r)}}=rp_E, \ \ q_{E^{(r)}}=rq_E, \end{align}
and if $0<p,q<\infty$, then \begin{align}\label{dis clu}(E^{(p)})^{(q)}=E^{(pq)}.\end{align}
If $E$ is a symmetric Banach function space and $p\geq1$, then $E^{(p)}$ is a symmetric Banach function space.

Let $E_i$ be a quasi-Banach function space for $i=1,2.$ The pointwise product space $E_1\odot E_2$ is defined by
$$E_1\odot E_2=\{f\in L_0(0,\infty): f=f_1f_2, f_i\in E_i, i=1,2\}$$
with functional $\|\cdot\|_{E_1\odot E_2}$ being defined by
$$\|f\|_{E_1\odot E_2}=\inf\{\|f\|_{E_1}\|f\|_{E_2}: f=f_1f_2, f_i\in E_i, i=1,2\}.$$

We need the following lemmas (see \cite[Theorem 1 (iii), Corollary 2]{ref8} and \cite[Theorem 6]{ref9}).

\begin{lemma}\label{lem2.1}
 Let $E$ and $F$ be two symmetric Banach function spaces.

{\rm (i)} If $0<p<\infty$, then $(E\odot F)^{(p)}=E^{(p)}\odot F^{(p)}.$

{\rm (ii)}   $L_1(0,\infty)=E\odot E^\times$.

{\rm (iii)}   If $1<p<\infty$, then $(E^{(p)})^{\times}=(E^\times)^{(p)}\odot L_{p'}(0,\infty).$
\end{lemma}

\begin{lemma}\label{lem2.2} Let $E$  be  a symmetric Banach function space  which is separable or has the Fatou property with $p_E>p$. Then $E^{(\frac{1}{p})}$ can be renormed as a symmetric Banach space.
\end{lemma}

\begin{proof}
 We first consider the case of $0<p\leq1$. Since $E$  is symmetric and  $\frac{1}{p}\geq1$, we have that $E^{(\frac{1}{p})}$ is a symmetric Banach function space (see \cite[p 53]{ref8}). For $p>1$, by \cite[Theorem 3.2]{ref4}, we have that $E$ is an interpolation space for the couple $\big(L_p(0,\infty), L_\infty(0,\infty)\big)$. It follows that $E^{(\frac{1}{p})}$ is an interpolation space for the couple $\big(L_1(0,\infty), L_\infty(0,\infty)\big)$ (see  \cite[Theorem 3.5]{ref4}). Thus according to  Lemma 2.2 in \cite{ref19}, we get that $E^{(\frac{1}{p})}$ can be renormed as a symmetric Banach function space.
\end{proof}

\subsection{ Noncommutative symmetric spaces and martingales}\label{subsetion3.2}
For a given noncommutative measure space $(\mathcal{M},\tau)$ and a symmetric quasi-Banach function space $(E,\|\cdot\|_E)$ on $(0,\infty)$, we define the
corresponding noncommutative symmetric space by setting
$$E(\mathcal{M},\tau)=\{x\in L_0(\mathcal{M}):\mu_t(x)\in E\}$$
equipped with the quasi-norm$$\|x\|_{E(\mathcal{M},\tau)}:=\|\mu_t(x)\|_E.$$
It is well-known that $E(\mathcal{M},\tau)$ (denoted by $E(\mathcal{M})$ for convenience) is a quasi-Banach space. Note that if $0< p<\infty$ and $E = L_p(0,\infty)$,  then
$E(\mathcal{M},\tau) = L_p(\mathcal{M},\tau)$ is the usual noncommutative $L_p$-space associated with $(\mathcal{M},\tau).$
We refer to \cite{ref1,ref19} for more details and historical references on these spaces.

\smallskip

Let $(\mathcal{M},\tau)$ be a noncommutative measure admitting a martingale structure, that is, there exists an increasing sequence of von Neumann
subalgebras $({\mathcal{M}}_n)_{n\geq1}$ of $\mathcal{M}$ such that the union of the ${\mathcal{M}}_n$'s is weak$^*$-dense in $\mathcal{M}$.
For each $n\geq1$, the unique
conditional expectation ${\mathcal{E}}_n$ from ${\mathcal{M}}$ onto ${\mathcal{M}}_n$ extends to a contractive projection from $L_p({\mathcal{M}},\tau)$ onto $L_p({\mathcal{M}}_n,\tau_n)$ for all $1\leq p\leq \infty$.  More
generally, if $E$ is a symmetric Banach function space on $(0,\infty)$ that belongs to $\mbox{Int}(L_1, L_\infty)$---the interpolation spaces, then ${\mathcal{E}}_n$ is bounded from $E({\mathcal{M}},\tau)$ onto $E({\mathcal{M}}_n,\tau_n)$.

A sequence $x=(x_n)_{n\geq1}$ is called a $E(\mathcal{M})$-martingale if $x_n\in E(\mathcal{M})$ and
${\mathcal{E}}_{n}(x_{n+1})=x_n$, $\forall n\geq1.$ In this case, we set $\|x\|_E=\sup_{n\geq1}\|x_n\|_E$. If $\|x\|_E<\infty$,  then $x$ is called a bounded  $E(\mathcal{M})$-martingale.
The column and row conditioned Hardy spaces $\mathsf{h}_E^c(\mathcal{M})$ and
  $\mathsf{h}_E^r(\mathcal{M})$ are defined to be respectively the completions of the space of all the finite $E(\mathcal{M})$-martingales under the associated quasi-norms
   $\|x\|_{\mathsf{h}_E^c}=\|s_c(x)\|_E$ and  $\|x\|_{\mathsf{h}_E^r}=\|s_r(x)\|_E$,
 where $s_c(x)$ and $s_r(x)$ are the column and row conditioned square functions of $x$, defined by (with the convention that $\mathcal{E}_0=\mathcal{E}_1$)
$$s_c(x)=(\sum\limits_{k\geq1}{\mathcal{E}}_{k-1}|dx_k|^2)^{1/2} \ \  \mbox{and} \ \  s_r(x)=s_c(x^*).$$

 In general, we have no explicit description of elements in $\mathsf{h}_E^c(\mathcal{M})$  or $\mathsf{h}_E^r(\mathcal{M})$, that is,  not all the elements of $\mathsf{h}_E^c(\mathcal{M})$ and $\mathsf{h}_E^r(\mathcal{M})$ can be represented by a martingale.  However, if $E=L_p(0,\infty)$ for $0<p<\infty$, then $\mathsf{h}_E^c(\mathcal{M})=\mathsf{h}_p^c(\mathcal{M})$ and $\mathsf{h}_E^r(\mathcal{M})=\mathsf{h}_p^r(\mathcal{M})$, namely the column and row conditioned Hardy spaces of noncommutative martingales.
  As remarked in \cite{ref15}, if $E$ is a symmetric Banach function space on $(0,\infty)$ with the Fatou property which is an interpolation space for the couple $\big(L_p(0,\infty), L_q(0,\infty)\big)$ for some $1<p<q<\infty$,  then every element of $\mathsf{h}_E^c(\mathcal{M})$ and $\mathsf{h}_E^r(\mathcal{M})$ can be represented by a martingale.

\subsection{Symmetric space moment inequality}
\label{subsection 3.3}
Starting with the $p$-moment inequalities---\eqref{JNlambdaMoment}, we are able to obtain the moment characterisation in terms of symmetric spaces in a more direct way. In the $\mathsf{bmo}^c(\mathcal M)$ case, our approach is much more efficient than the interpolation arguments given in \cite{ref19}.

 We first introduce the symmetric Lipschitz spaces $\Lambda^c_{\beta,E}(\mathcal{M})$ and $\Lambda^r_{\beta,E}(\mathcal{M})$.

 \begin{definition} \label{lem4.1}Let $\beta\geq0$. Let $E$ be  a separable  symmetric Banach  function space on $(0,\infty)$ with
$1<p_E\leq{q}_E<\infty$. We define
$$\Lambda^c_{\beta,E}(\mathcal{M})=\{x\in L_2(\mathcal{M}): \|x\|_{\Lambda^c_{\beta,E}}<\infty\},$$ where
$$
\|x\|_{\Lambda^c_{\beta,E}}=\max\{\|\mathcal{E}_1(x)\|_\infty,\ \ \sup_n
\sup_{e\in{\mathcal{P}(\mathcal{M}_n)}}{(\tau(e))}^{-\beta}\|(x-x_n)\frac{e}{\|e\|_E}\|_{\mathsf{h}_E^c}\}
$$
and $$ \Lambda^r_{\beta,E}(\mathcal{M})=\{x\in L_2(\mathcal{M}): x^*\in \Lambda^c_{\beta,E}(\mathcal{M})\}.$$
 Let $ \Lambda^{0,c}_{\beta,E}(\mathcal{M})$  and $ \Lambda^{0,r}_{\beta,E}(\mathcal{M})$ be their subspaces of all $x$ with $\mathcal{E}_1(x)=0$.
\end{definition}

\begin{theorem}  \label{lem4.2}Let $\beta\geq0$ and $x\in \Lambda^c_{\beta}(\mathcal{M})$.  Let $E$ be  a separable  symmetric Banach function space on $(0,\infty)$ with
$1<p_E\leq{q}_E<\infty$. Then we have that
$$\alpha_E^{-1}\|x\|_{{\Lambda}^c_{\beta}}\leq\|x\|_{{\Lambda^c_{\beta,E}}}\leq\beta_E\|x\|_{{\Lambda}^c_{\beta}}.$$
The constants $\alpha_E$ and $\beta_E$ depend on $E$.
The same inequalities hold for $\|\cdot\|_{{\Lambda}^r_{\beta}}$ and $\|\cdot\|_{{\Lambda}^r_{\beta,E}}$.
\end{theorem}

\begin{lemma} \label{lem4.3}  Let $E$ be a separable symmetric Banach function space on $(0,\infty)$  with $1<p_E\leq{q}_E<\infty$.
Choose $p$ such that ${q}_E<p$. Then we have that $E=F^{(p')}\odot L_p(0,\infty)$, where
$F=(E^{\times(\frac{1}{p'})})^{\times}$. Here and below, for $1\leq p\leq\infty$, $p'$ denotes the conjugate index of $p$.
\end{lemma}

\begin{proof}  By \eqref{dis cll}, it follows that $1<p'<p_{E^\times}\leq{q}_{E^\times}<\infty$. Since $E^{\times}$ is  symmetric and  has the Fatou property,  by Lemma 2.2,  $E^{\times(\frac{1}{p'})}$ can be renormed as a  symmetric Banach function space.
Then  by \eqref{dis clu} and (iii) of  Lemma \ref {lem2.1}, we have that
$$E=(E^\times)^\times=([E^{\times(\frac{1}{p'})}]^{(p')})^\times=([E^{\times(\frac{1}{p'})}]^{\times})^{(p')}\odot
L_p(0,\infty)=F^{(p')}\odot L_p(0,\infty).$$
\end{proof}

We will use repeatedly the following fact which follows from Theorem 2 of \cite{ref10}.

\begin{lemma} \label{lem4.4} Let $E$ and $F$ be two symmetric Banach function spaces on $(0,\infty)$. Then $E\odot F$ is a symmetric quasi-Banach function space on $(0,\infty)$ and the following formula holds:
$$\|\chi_{[0,t]}\|_{E\odot F}=\|\chi_{[0,t]}\|_{E}\|\chi_{[0,t]}\|_{ F} \ \   \mbox{for} \ \ t\in (0,\infty).$$
\end{lemma}

\begin{proof}[Proof of Theorem 3.4.]  We may assume $\mathcal E_1(x)=0$.
Let $F$, $p$ and $p'$ be as in Lemma \ref{lem4.3}. Fix $n$ and $P\in{\mathcal{P}(\mathcal{M}_n)}$. By definition $\mu_t(P)=\chi_{[0,\tau(P)]}(t)$, and thus by
Lemma \ref{lem4.3} and  Lemma \ref{lem4.4}, we have
$$\|P\|_E=\|P\|_{F^{(p')}}\|P\|_p.$$
Therefore,  by the equality $E^{\frac{1}{2}}=F^{(\frac{p'}{2})}\odot L_{\frac{p}{2}}(0,\infty)$,
 we have that \begin{eqnarray}\label{dis cpp}
(\tau(P))^{-\beta}\|(x-x_n)\frac{P}{\|P\|_E}\|_{\mathsf{h}_E^c}
&=&(\tau(P))^{-\beta}\big\|\frac{P}{\|P\|_E}{s}_c^2(x-x_n)\frac{P}{\|P\|_E}\big\|_{E^{(\frac{1}{2})}}^{\frac{1}{2}}\nonumber\\
&\leq&(\tau(P))^{-\beta}\big\|\frac{P}{\|P\|_E^2}\big\|_{F^{(\frac{p'}{2})}}^{\frac{1}{2}}\|P
{s}_c^2(x-x_n)P\|_{\frac{p}{2}}^{\frac{1}{2}}\nonumber\\
&=&(\tau(P))^{-\beta}\frac{1}{(\tau(P))^{\frac{1}{p}}}\|(x-x_n)P\|_{\mathsf{h}_p^c}
\end{eqnarray}
 which implies
that $\|x\|_{{\Lambda^c_{\beta,E}}}\leq \|x\|_{{\Lambda}^c_{\beta,p}}.$ Then  by (iii) of Theorem \ref{thm:JNlambda}, we have that
$\|x\|_{{\Lambda^c_{\beta,E}}}\leq \beta_E\|x\|_{{\Lambda}^c_{\beta}}.$

Conversely, choose $p$ such that $p<p_E$. Then by Lemma \ref{lem2.2},  $E^{(\frac{1}{p})}$ can be renormed as a symmetric Banach function space.  Thus by (ii) of Lemma \ref{lem2.1}, we have that $L_1(0,\infty)=E^{(\frac{1}{p})}\odot E^{(\frac{1}{p})\times}$.
It follows that
$L_p(0,\infty)=E\odot F$ by (i) of Lemma \ref{lem2.1},
where $F=(E^{(\frac{1}{p})\times})^{(p)}$. Hence, for all ${e}\in \mathcal{P}(\mathcal{M}_n)$,  by Lemma \ref{lem4.4}, we have
\begin{eqnarray}\label{dis cpw}
(\tau(e))^{\frac{1}{p}}=\|e\|_E\|e\|_F.
\end{eqnarray}
  Using the equality
$L_{\frac{p}{2}}(0,\infty)=E^{(\frac{1}{2})}\odot F^{(\frac{1}{2})}$, one has
\begin{eqnarray}\label{dis cpq}
\|(x-x_n)\frac{e}{(\tau(e))^{\frac{1}{p}}}\|_{\mathsf{h}_p^c}
&=&\big\|\frac{e}{\|e\|_F}\frac{e}{\|e\|_E}{s}_c^2(x-x_n)\frac{e}{\|e\|_E}\frac{e}{\|e\|_F}\big\|_{\frac{p}{2}}^{\frac{1}{2}}\nonumber\\
&\leq&\big\|\frac{e}{\|e\|_F^2}\big\|_{F^{(\frac{1}{2})}}^{\frac{1}{2}}\big\|\frac{e}{\|e\|_E}{s}_c^2(x-x_n)\frac{e}{\|e\|_E}\big\|_{E^{(\frac{1}{2})}}^{\frac{1}{2}}\nonumber\\
&=&\big\|(x-x_n)\frac{e}{\|e\|_E}\big\|_{\mathsf{h}_E^c}
\end{eqnarray}
which implies $\|x\|_{{\Lambda^c_{\beta,p}}}\leq\|x\|_{{\Lambda}^c_{\beta,E}} $.
Therefore, by  Theorem \ref{thm:JNlambda},  we get the desired result
$$\|x\|_{{\Lambda^c_{\beta}}}\leq \alpha_E\|x\|_{{\Lambda}^c_{\beta,E}}.$$
\end{proof}

\begin{remark}\label{lem4.5}
When $0<p_E\leq{q}_E\leq 1$, it is easy to obtain the left inequality  from the proof of Theorem \ref{lem4.2}.
\end{remark}

\section{The  atomic Hardy space $\mathsf{h}_{p,E}(\mathcal{M})$ for $0<p\leq1$}
\label{section 4}
In this section, we show that the noncommutative Hardy space $\mathsf{h}^c_{p,E}(\mathcal{M})$ ($0<p\leq1$) defined via symmetric space atoms is also a predual space of Lipshitz space ${\Lambda}^c_{\frac1p-1}(\mathcal{M})$.

\begin{definition} \label{lem4.6} Let $0<p\leq1$. Let $E$ be  a  separable symmetric Banach function space on $(0,\infty)$ and
$1<p_E\leq{q}_E<\infty$.  An element $a\in L_p(\mathcal{M})$ is called a $(p,E)_c$-atom with respect to $(\mathcal{M}_n)_{n\geq1}$, if there exist $n\geq1$ and a projection $e\in \mathcal{M}_n$ such that:

(i) $\mathcal{\mathcal{E}}_n(a)=0;$

(ii) $r(a)\leq e;$

(iii) $\|a\|_{\mathsf{h}_E^c}\leq (\tau(e))^{1-1/p}\|e\|_{E^\times}^{-1}.$

Replacing (ii) by $l(a)\leq e$, we have the notion of a $(p,E)_{r}$-atom.
\end{definition}
Note that if $E=L_q(0,\infty)$ for  $1<q< \infty$, they are  $(p,q)$-atoms defined in \cite[Definition 4.1]{ref3}.

\begin{definition}  \label{lem4.8} Let $0<p\leq1$.  Let $E$ be  a  separable symmetric Banach function space on $(0,\infty)$  and
$1<p_E\leq{q}_E<\infty$.   We define $\mathsf{h}_{p,E}^c(\mathcal{M})$  as the space of all operators $x\in L_p(\mathcal{M})$ which admits a decomposition
$$x=\sum_k\lambda_ka_k,$$
where for each $k$, ${a_k}$ is either a $(p,E)_c$-atom or an element of the unit ball of $L_p(\mathcal{M}_1)$, and $\lambda_k\in \mathbb{C}$ satisfying $\sum_k|\lambda_k|^p<\infty$.  We equip this space with the $p$-norm:
$$\|x\|_{\mathsf{h}_{p,E}^c}=\inf{(\sum_k|\lambda_k|^p)^{\frac{1}{p}}},$$
where the infimum is taken over all decompositions of $x$ described above.
We also define  the subspace:
 $$\mathsf{h}_{p,E}^{0,c}(\mathcal{M})=\{x\in \mathsf{h}_{p,E}^c(\mathcal{M}): \  \mathcal{E}_1(x)=0\}.$$
Similarly, we define $\mathsf{h}_{p,E}^r(\mathcal{M})$ and $\mathsf{h}_{p,E}^{0,r}(\mathcal{M})$.
\end{definition}

\begin{lemma}\label{lem4.7} Let $0<p\leq1$.  If $a$ is a $(p,E)_{c}$-atom, then
$$\|a\|_{\mathsf{h}_p^c}\leq 1.$$
The similar inequality holds for $(p,E)_{r}$-atom.
\end{lemma}
\begin{proof} Let $e$ be a projection associated with $a$ satisfying (i)-(iii) of Definition \ref{lem4.6}.
Then $s_c(a)=es_c(a)=s_c(a)e$ (see \cite[Proposition2.2]{ref2}). Nothing that $p<p_E$, similar with the proof of \eqref{dis cpq},
\begin{eqnarray*}
\|a\|_{\mathsf{h}_p^c}
&\leq&\|s_c(a)\|_{E^{(\frac{1}{2})}}^{\frac{1}{2}}\|e\|_F\\
&=&\|a\|_{\mathsf{h}_E^c}\frac{(\tau(e))^{\frac{1}{p}}}{\|e\|_E} \  (\text{by}\ \  \eqref{dis cpw})\\
&\leq&\frac{\tau(e)}{\|e\|_{E^\times}\|e\|_E}  \  \text{ (by (iii) of Definition  \ref{lem4.6} })\\
&=&1 \ \  \text{ ( by $L_1(0,\infty)=E\odot E^\times$ and Lemma  \ref{lem4.4} )}.
\end{eqnarray*}
Thus we obtain the desired result.
\end{proof}

\begin{lemma}  \label{lem4.10} Let $E$ be  a  separable symmetric Banach function space on $(0,\infty)$ and $1<p_E\leq{q}_E<\infty$. For a given $q_E<q_0<\infty$ and $q=\max\{2,q_0\}$,  $L_q(\mathcal{M})$  embeds densely and continuously into $\mathsf{h}_{p,E}^c(\mathcal{M})$.
\end{lemma}
\begin{proof} Let $x\in L_q(\mathcal{M})$. We decompose it as a linear combination of two atoms:
$$x=C_qC_E\|x-\mathcal{E}_1(x)\|_q\frac{x-\mathcal{E}_1(x)}{C_qC_E\|x-\mathcal{E}_1(x)\|_q}+\|\mathcal{E}_1(x)\|_q\frac{\mathcal{E}_1(x)}
{\|\mathcal{E}_1(x)\|_q},$$
where $C_E$ is the constant in the inequality $\|y\|_E\leq C_E\|y\|_q$ for all $y\in L_q(\mathcal{M})$, and $C_q$ is the constant in the noncommutative Burkholder inequality $\|y\|_{\mathsf{h}_q^c}\leq C_q\|y\|_q$ for all $y\in L_q(\mathcal{M})$ (see \cite{ref7}). Then
${\mathcal{E}_1(x)}/{\|\mathcal{E}_1(x)\|_q}\in L_q(\mathcal{M}_1)\subset L_1(\mathcal{M}_1)$ and
$$\big\|\frac{\mathcal{E}_1(x)}{\|\mathcal{E}_1(x)\|_q}\big\|_1=\frac{\|{\mathcal{E}_1(x)}\|_1}{{\|\mathcal{E}_1(x)\|_q}}\leq1.$$
 Also,
$$\frac{x-\mathcal{E}_1(x)}{C_qC_E\|x-\mathcal{E}_1(x)\|_q}=\frac{x-\mathcal{E}_1(x)}{C_qC_E\|x-\mathcal{E}_1(x)\|_q}\cdot 1\doteq ae.$$
Clearly, $\mathcal{E}_1(a)=0$ and $$\big\|\frac{x-\mathcal{E}_1(x)}{C_qC_E\|x-\mathcal{E}_1(x)\|_q}\big\|_{\mathsf{h}_E^c}\leq
\frac{C_E\|x-\mathcal{E}_1(x)\|_{\mathsf{h}_q^c}}{C_qC_E\|x-\mathcal{E}_1(x)\|_q}\leq1.$$
Thus $$\|x\|_{\mathsf{h}_{p,E}^c}\leq C_qC_E\|x-\mathcal{E}_1(x)\|_q+\|\mathcal{E}_1(x)\|_q\leq (2C_qC_E+1)\|x\|_q.$$
The density is trivial.
\end{proof}

 \begin{proposition}  \label{lem4.11} Let $0<p\leq1$ and $\beta=1/p-1$. Let $E$ be  a separable  symmetric Banach function space on $(0,\infty)$. If $1<p_E\leq{q}_E<\infty$, then
$\big(\mathsf{h}_{p,E}^{0,c}(\mathcal{M})\big)^*=\Lambda_{\beta,E^\times}^{0,c}(\mathcal{M})$ with equivalent norms. More precisely,

{\rm(i)} Every $y\in \Lambda_{\beta,E^\times}^{0,c}(\mathcal{M})$ defines a continuous linear functional on $\mathsf{h}_{p,E}^{0,c}(\mathcal{M})$ by
\begin{eqnarray}\label{dis cpa}
\varphi_y(x)=\tau(y^*x)
\end{eqnarray}
 for all $x\in \mathsf{h}_{p,E}^{0,c}(\mathcal{M})$, and $\|\varphi_y\|_{\big(\mathsf{h}_{p,E}^{0,c}(\mathcal{M})\big)^*}\leq \|y\|_{ \Lambda_{\beta,E^\times}^{0,c}}$;

{\rm(ii)} Conversely, each $\varphi\in \big(\mathsf{h}_{p,E}^{0,c}(\mathcal{M})\big)^*$ is given as \eqref{dis cpa} by some $y\in \Lambda_{\beta,E^\times}^{0,c}(\mathcal{M})$, and $\|y\|_{ \Lambda_{\beta,E^\times}^{0,c}}\leq C_E\|\varphi_y\|_{\big(\mathsf{h}_{p,E}^{0,c}(\mathcal{M})\big)^*}$.

Similarly, $\big(\mathsf{h}_{p,E}^{0,r}(\mathcal{M})\big)^*=\Lambda_{\beta,E^\times}^{0,r}(\mathcal{M})$ with equivalent norms.
 \end{proposition}

\begin{proof} (i) Let $y\in \Lambda_{\beta,E^\times}^{0,c}(\mathcal{M})$. If $a$ is  a  $(p,E)_c$-atom with $\mathcal{\mathcal{E}}_n(a)=0$ for some
$n\geq1$ and $a=ae$ for some projection $e\in \mathcal{M}_n$ satisfying $\|a\|_{\mathsf{h}_E^c}\leq (\tau(e))^{1-1/p}\|e\|_{E^\times}^{-1},$ then using the duality inclusion
$\mathsf{h}^c_{E^\times}(\mathcal{M})\subset\big(\mathsf{h}^c_E(
\mathcal{M})\big)^*$ with constant 1 (see \cite[Theorem 2.5]{ref19}), we have that
\begin{eqnarray*}
|\tau(y^*a)|
&=&|\tau\big((y-y_n)^*ae\big)|\\
&\leq& \|(y-y_n)e\|_{\mathsf{h}^c_{E^\times}}\|a\|_{\mathsf{h}_E^c}\\
&\leq&\|(y-y_n)e\|_{\mathsf{h}^c_{E^\times}} (\tau(e))^{1-1/p}\|e\|_{E^\times}^{-1}\\
&\leq&\|y\|_{\Lambda_{\beta,E^\times}^{0,c}}.
\end{eqnarray*}

(ii) Let $\varphi\in \big(\mathsf{h}_{p,E}^{0,c}(\mathcal{M})\big)^*$. Set $q=\max\{q_0,2\},$ where $q_E<q_0<\infty$. By Lemma  \ref{lem4.10}, we can find $y\in L_{q'}(\mathcal{M})$  such that
$$\varphi(x)=\tau(y^*x), \ \  \forall x\in L_{q}(\mathcal{M}).$$
Fix $n\geq1$ and $e\in \mathcal{P}(\mathcal{M}_n)$. For a fixed arbitrary and small enough $\varepsilon>0$, by the inclusion $\big(\mathsf{h}^c_E(
\mathcal{M})\big)^*\subset\mathsf{h}^c_{E^\times}(\mathcal{M})$ with contant $C_E$, we may choose $x\in L_{q}(\mathcal{M})$ such that $\|x\|_{\mathsf{h}_{E}^c}\leq 1$ so that
$$C_E|\tau(e(y-y_n)^*x)|+\varepsilon\geq\|(y-y_n)e\|_{\mathsf{h}_{E^\times}^c}.$$
Clearly, we may assume that $\mathcal{E}_n(x)=0$ and $xe=x$. Set
$$a=\frac{x}{\|x\|_{\mathsf{h}_{E}^c}(\tau(e))^{1/p-1}\|e\|_{E^\times}}.$$
Then $a$ is a $(p,E)_c$-atom and
\begin{eqnarray*}
\|\varphi\|
&\geq&|\tau\big((y-y_n)^*a\big)|\\
&=& \frac{1}{\|x\|_{\mathsf{h}_{E}^c}(\tau(e))^{1/p-1}\|e\|_{E^\times}}|\tau(e(y-y_n)^*x)|\\
&\geq&\frac{1}{C_E(\tau(e))^\beta\|e\|_{E^\times}}(\|(y-y_n)e\|_{\mathsf{h}_{E^\times}^c}-\varepsilon).
\end{eqnarray*}
By the arbitrariness of $\varepsilon$, and taking supremum over $n$ and $e\in \mathcal{P}(\mathcal M_n)$, we get $C_E\|\varphi\|\geq \|y\|_{\Lambda_{\beta,E^\times}^{0,c}}.$
\end{proof}

By Proposition \ref{lem4.11} and Theorem \ref{lem4.2}, we arrive at the main result of this subsection.

\begin{theorem}\label{thm:predual}
Let $0<p\leq1$ and $\beta=1/p-1$. Let $E$ be  a separable  symmetric Banach function space on $(0,\infty)$. If $1<p_E\leq{q}_E<\infty$, then
$\big(\mathrm{h}_{p,E}^{0,c}(\mathcal{M})\big)^*=\Lambda_{\beta}^{0,c}(\mathcal{M})$ with equivalent norms. The similar results hold for the row and mixture spaces.
\end{theorem}

\begin{remark}
 By Lemma  \ref{lem4.7}, we have the obvious inclusion $\mathsf{h}_{p,E}^c(\mathcal{M})\subset \mathsf{h}_p^c(\mathcal{M})$ for $0<p\leq1$. When $p=1$, the converse inclusion is also true \cite{ref19}; it is, however, an open question in the cases $0<p<1$ even though they have the same dual from Theorem \ref{thm:predual} and \cite[Theorem 5.4]{ref3}.
 \end{remark}

\begin{remark} \label{lem4.13}
 We may consider the crude symmetric atoms. Let $0<p\leq1$. Let $E$ be  a  separable  symmetric Banach function space on $(0,\infty)$ with
$1<p_E\leq{q}_E<\infty$.  An element $a\in L_p(\mathcal{M})$ is called a $(p,E)_c$-crude atom with respect to $(\mathcal{M}_n)_{n\geq1}$, if there exist $n\geq1$ and a factorization $a=yb$ such that:

(i) $\mathcal{\mathcal{E}}_n(y)=0;$

(ii) $b\in L_{\frac{p}{1-p}}(\mathcal{M}_n), \ b\in E^\times(\mathcal{M}_n) \ \mbox{and} \  \|b\|_{{\frac{p}{1-p}}}\|b\|_{E^\times}\leq1 $ for $0<p<1$;\ \  $ b\in E^\times(\mathcal{M}_n) \ \mbox{and} \  \|b\|_{E^\times}\leq1 $ for $p=1$;

(iii) $\|y\|_{\mathsf{h}_E^c}\leq 1.$

\noindent{Similarly, we define the notion of a $(p,E)_{r}$-crude atom with $a=yb$ replaced by $a=by$.}

 Similar to Definition  \ref{lem4.8}, we define $\mathsf{h}_{p,E;\rm{crude}}^c(\mathcal{M})$  based on the $(p,E)_{c}$-crude atoms as
building blocks.

Note that if $a$ is a $(p,E)_{c}$-atom with the associated projection $e\in \mathcal{M}_n$, then $a$ is a $(p,E)_{c}$-crude atom. Indeed,
for $0<p<1$, we write $$a=(\|e\|_{{\frac{p}{1-p}}}\|e\|_{E^\times}a)\frac{e}{\|e\|_{{\frac{p}{1-p}}}\|e\|_{E^\times}}=yb;$$
for $p=1$,
we write $a=(\|e\|_{E^\times}a)\frac{e}{\|e\|_{E^\times}}=yb$.

One can check that $\mathsf{h}_{p,E;\rm{crude}}^c(\mathcal{M})\subset \mathsf{h}_p^c(\mathcal{M})$ for $0<p\leq1$. So we have that
$$\mathsf{h}_{p,E}^c(\mathcal{M})\subset \mathsf{h}_{p,E;\rm{crude}}^c(\mathcal{M})\subset \mathsf{h}_p^c(\mathcal{M})$$ for $0<p\leq1$.
As commented the previous remark, it is an interesting question to show the three spaces are the same. We shall take care of it elsewhere.
\end{remark}

\bibliographystyle{amsalpha}

\begin{thebibliography}{99}



\bibitem{ref1}T. Bekjan, Duality for symmetric Hardy spaces of noncommutative martingales, Math. Zeit. 289 (2018), 787-802.
 \bibitem{ref2} T. Bekjan, Z. Chen, M. Perrin, Z. Yin, Atomic decomposition and interpolation for Hardy spaces of noncommutative martingales, J. Funct. Anal. 258 (7) (2010),  2483-2505.
   \bibitem{ref19} T. Bekjan, Z. Chen,  M. Raikhan, M. Sun,  Interpolation and John-Nirenberg inequality on symmetric spaces of noncommutative martingales,  Studia Math, (2021). DOI: 10.4064/sm200508-11-12.
\bibitem{ref3}Z. Chen, N. Randrianantoanina, Q. Xu, Atomic decompositions for noncommutative martingales, arXiv: 2001.08775v1, math.
    OA(2020).
    \bibitem{ref39} I. Cuculescu, Martingales on von Neumann algebras. Journal of Multivariate Analysis,  1(1)(1971), 17-27.
\bibitem{ref4}S. Dirksen, B. dePagter, D. Potapov, and F. Sukochev, Rosenthal inequalities in noncommutative symmetric spaces, J. Funct. Anal. 261 (2011), 2890-2925.
    \bibitem{ref20} T. Holmstedt,  Interpolation of quasi-normed spaces, Math. Scand. 26 (1970), 177-199.
\bibitem{ref5}G. Hong, T. Mei, John-Nirenberg inequality and atomic decomposition for noncommutative martingales, J. Funct. Anal.
    263 (2012), 1064-1097.
\bibitem{JRWZ20} Y. Jiao, N. Randrianantoanina, L. Wu, D. Zhou,  Square Functions for Noncommutative Differentially Subordinate Martingales, Commun. Math. Phys. 374  (2020), 975-1019.
  \bibitem{JOW18}  Y. Jiao, A. Os{\c e}kowski, L. Wu, Inequalities for noncommutative differentially subordinate martingales, Adv. Math. 337 (2018), 216-259.
\bibitem{ref6}Y. Jiao, A. Os{\c e}kowski, L. Wu, Noncommutative good-$\lambda$ ineaualities, arXiv: 1802.07057v1, math.
    OA(2018).
  \bibitem{ref30}  M. Junge, Doob's inequality for noncommutative martingales, J. Reine Angew. Math. 549 (2002), 149-190.
  \bibitem{ref31} M. Junge, Q. Xu, Noncommutative Burkholder/Rosenthal inequalities, Ann. Probab. 31 (2) (2003), 948-995.
   \bibitem{ref32} M. Junge, Q. Xu, Noncommutative maximal ergodic theorems, J. Amer. Math. Soc. 20 (2006), 385-439.
     \bibitem{ref7}M. Junge, Q. Xu, Noncommutative Burkholder/Rosenthal inequalities II: Applications, IsraelJ. Math. 167 (2008), 227-282.
    \bibitem{ref18} M. Junge, M. Musat, A noncommutative version  of John-Nirenberg theorem, Trans. Amer. Math. Soc. 359(1) (2007), 115-142.
\bibitem{ref8}P. Kolwicz, K. Le\' snik, L. Maligranda, Pointwise products of some Banach function spaces and
    factorization, J. Funct. Anal. 266 (2014), 616-659.

  \bibitem{ref9}  G. Y. Lozanovskii, On some Banach lattices, Siberian Math. J. 10 (1969), 419-430.
   \bibitem{ref10} J. Lindenstrauss, L. Tzafriri, Classical Banach Spaces. II: Function Spaces, Ergeb. Math. Grenzgeb.
  (Results in Mathematics and Related Areas), vol.97, Springer-Verlag, Berlin, 1979.
  \bibitem{ref11}  M. Musat, Interpolation between noncommutative BMO and noncommutative $L_p$-spaces, J. Funct. Anal.
      202 (2003), 195-225.
 \bibitem{ref36} J. Parcet, N. Randrianantoanina, Gundy's decomposition for noncommutative martingales and applications, Proc. London Math. Soc. 93 (2006), 227-252.
\bibitem{ref12} J. Peetre,   G. Sparr, Interpolation and non-commutative integration, Ann. Mat. Pura. Appl. Vol. CIV (1975), 187-207.
\bibitem{ref38} M. Perrin, A noncommutative Davis' decomposition for martingales, J. London Math. Soc. 80 (3) (2009), 627-648.
\bibitem{ref35} G. Pisier, Q. Xu, Noncommutative martingale inequalities, Comm. Math. Phys. 189 (1997), 667-698.
\bibitem{ref13}N. Randrianantoanina, Conditioned square functions for noncommutative martingales, Ann. Probab. 35 (2007), 1039-1070.
\bibitem{ref33} N. Randrianantoanina, Noncommutative martingale transform, J. Funct. Anal. 194 (1) (2002), 181-212.
\bibitem{ref14} N. Randrianantoanina, L. Wu, Martingale inequalities in noncommutative symmetric spaces, J. Funct. Anal.
    269 (2015), 2222-2253.

 \bibitem{ref15}   N. Randrianantoanina, L. Wu, Q. Xu, Noncommutative Davis type decompositions and applications, J. London Math. Soc. 99 (2019), 97-126.
\bibitem{ref25}F. Weisz,  Martingale Hardy Spaces and Their Applications in Fourier Analysis. Berlin: Springer-Verlang, 1994.
    \bibitem{ref16}Q. Xu, Analytic functions with values in lattices and symmetric spaces of measurable operators, Math. Proc. Camb.
        Phil. Soc. 109 (1991), 541-563.

 \bibitem{ref17} Q. Xu, Applications du th$\mbox{\'{e}}$or$\mbox{\`{e}}$ me de factorization pour les fonctions $\mbox{\`{a}}$ valeurs op$\mbox{\'{e}}$rateurs, Studia Math. 95 (1990), 273-292.
    \bibitem{ref26}  Q. Xu, Noncommtative $L_p$-Spaces and Martingale Inequalities, book manuscript, (2007).



\end{thebibliography}

\end{document}